# Lie–Butcher series, Geometry, Algebra and Computation

Hans Z. Munthe-Kaas and Kristoffer K. Føllesdal

**Abstract** Lie–Butcher (LB) series are formal power series expressed in terms of trees and forests. On the geometric side LB-series generalizes classical B-series from Euclidean spaces to Lie groups and homogeneous manifolds. On the algebraic side, B-series are based on pre-Lie algebras and the Butcher-Connes-Kreimer Hopf algebra. The LB-series are instead based on post-Lie algebras and their enveloping algebras. Over the last decade the algebraic theory of LB-series has matured. The purpose of this paper is twofold. First, we aim at presenting the algebraic structures underlying LB series in a concise and self contained manner. Secondly, we review a number of algebraic operations on LB-series found in the literature, and reformulate these as recursive formulae. This is part of an ongoing effort to create an extensive software library for computations in LB-series and B-series in the programming language Haskell.

## 1 Introduction

Classical B-series are formal power series expressed in terms of rooted trees (connected graphs without any cycle and a designated node called the root). The theory has its origins back to the work of Arthur Cayley [5] in the 1850s, where he realized that trees could be used to encode information about differential operators. Being forgotten for a century, the theory was revived through the efforts of understanding numerical integration algorithms by John Butcher in the 1960s and '70s [2, 3]. Ernst Hairer and Gerhard Wanner [15] coined the term *B-series* for an infinite formal series of the form

H. Z. Munthe-Kaas, K. Føllesdal
Department of Mathematics, University of Bergen, P.O. Box 7803, N-5020 Bergen, e-mail: hans.munthe-kaas@uib.no, e-mail: kristoffer.follesdal@uib.no

$$B_f(\alpha,y,t) := y + \sum_{\tau\in T} \frac{t^{|\tau|}}{\sigma(\tau)}\langle a,\tau\rangle \mathscr{F}_f(\tau)(y),$$

where $y \in \mathbb{R}^n$ is a 'base' point, $f\colon \mathbb{R}^n \to \mathbb{R}^n$ is a given vector field, $T = \{$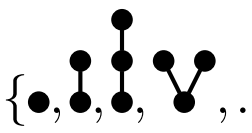$,\ldots\}$ is the set of rooted trees, $|\tau|$ is the number of nodes in the tree, $\alpha\colon T \to \mathbb{R}$ are the coefficients of a given series and $\langle \alpha, \tau\rangle \in \mathbb{R}$ denotes evaluation of $\alpha$ at $\tau$. The bracket hints that we later want to consider $\langle \alpha, \cdot\rangle$ as a linear functional on the vector space spanned by $T$. The animal $\mathscr{F}_f(\tau)\colon \mathbb{R}^n \to \mathbb{R}^n$ denotes special vector fields, called *elementary differentials*, which can be expressed in terms of partial derivatives of $f$. The coefficient $\sigma(\tau) \in \mathbb{N}$ is counting the number of symmetries in a given tree. This symmetry factor could have been subsumed into $\alpha$, but is explicitly taken into the series due to the underlying algebraic structures, where this factor comes naturally. The B-series $t \mapsto B_f(\alpha,y,t)$ can be interpreted as a curve starting in $y$. By choosing different functions $\alpha$, one may encode both the analytical solution of a differential equation $y'(t) = f(y(t))$ and also various numerical approximations of the solution.

During the 1980s and 1990s B-series evolved into an indispensable tool in analysis of numerical integration for differential equations evolving on $\mathbb{R}^n$. In the mid-1990s interest rose in the construction of numerical integration on Lie groups and manifolds [18, 16], and from this a need to interpret B-series type expansions in a differential geometric context, giving birth to *Lie–Butcher series* (LB-series), which combines B-series with Lie series on manifolds. It is necessary to make some modifications to the definition of the series to be interpreted geometrically on manifolds:

- We cannot add a point and a tangent vector as in $y + \mathscr{F}_f(\tau)$. Furthermore, it turns out to be very useful to regard the series as a Taylor-type series for the mapping $f \mapsto B_f$, rather than a series development of a curve $t \mapsto B_f(a,y,t)$. The target space of $f \mapsto B_f$ is differential operators, and we can remove explicit reference to the base point $y$ from the series.
- The mapping $f \mapsto B_f$ inputs a vector field and outputs a series which may represent either a vector field or a solution operator (flow map). Flow maps are expressed as a series in higher order differential operators. We will see that trees encode first order differential operators. Higher order differential operators are encoded by products of trees, called *forests*. We want to also consider series in linear combinations of forests.
- We will in the sequel see that the elementary differential map $\tau \mapsto \mathscr{F}_f(\tau)$ is a *universal arrow* in a particular type of algebras. The existence of such a uniquely defined map expresses the fact that the vector space spanned by trees (with certain algebraic operations) is a universal object in this category of algebras. Thus the trees encode faithfully the given algebraic structure. We will see that the algebra comes naturally from the geometric properties of a given *connection* (covariant derivation) on the manifold. For Lie groups the algebra of the natural connection is encoded by *ordered* rooted trees, where the ordering of the branches is important. The ordering is related to a non-vanishing torsion of the connection.

- The symmetry factor $\sigma(\tau)$ in the classical B-series is related to the fact that several different ordered trees correspond to the same unordered tree. This factor is absent in the Lie–Butcher series.
- The time parameter $t$ is not essential for the algebraic properties of the series. Since $\mathscr{F}_{tf}(\tau) = t^{|\tau|}\mathscr{F}_f(\tau)$, we can recover the time factor through the substitution $f \mapsto tf$.

We arrive at the definition of an abstract Lie–Butcher series simply as

$$\sum_{\omega \in \mathrm{OF}} \langle \alpha, \omega \rangle \omega, \tag{1}$$

where

OF = {$\mathbb{I}$, 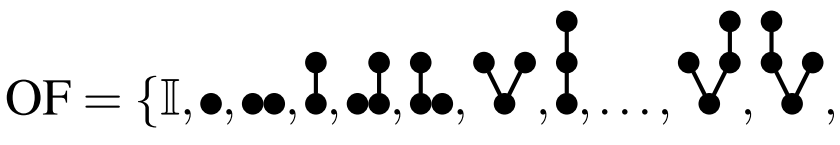 , ...}

denotes the set of all ordered forests of ordered trees, $\mathbb{I}$ is the empty forest, and $\alpha\colon \mathrm{OF} \to \mathbb{R}$ are the coefficients of the series. This abstract series can be mapped down to a concrete algebra (e.g. an algebra of differential operators on a manifold) by a universal mapping $\omega \mapsto \mathscr{F}_f(\omega)$.

We can identify the function $\alpha\colon \mathrm{OF} \to \mathbb{R}$ with its series (1) and say that a Lie–Butcher series $\alpha$ is an element of the graded dual vector space of the free vector space spanned by the forests of ordered rooted trees. However, to make sense of this statement, we have to attach algebraic and geometric meaning to the vector space of ordered forests. This is precisely explained in the sequel, where we see that the fundamental algebraic structures of this space arise because it is the universal enveloping algebra of a free post-Lie algebra. Hence we arrive at the precise definition:

*An abstract Lie–Butcher series is an element of the dual of the enveloping algebra of the free post-Lie algebra.*

We will in this paper present the basic geometric and algebraic structures behind LB-series in a self contained manner. Furthermore, an important goal for this work is to prepare a software package for computations on these structures. For this purpose we have chosen to present all the algebraic operations by recursive formulae, ideally suited for implementation in a functional programming language. We are in the process of implementing this package in the Haskell programming language. The implementation is still at a quite early stage, so a detailed presentation of the implementation will be reported later.

## 2 Geometry of Lie–Butcher series

B-series and LB-series can both be viewed as series expansions in a connection on a fibre bundle, where B-series are derived from the canonical (flat and torsion free) connection on $\mathbb{R}^n$ and LB-series from a flat connection with constant torsion on a fibre bundle. Rather than pursuing this idea in an abstract general form, we will provide insight through the discussion of concrete and important examples.

## 2.1 Parallel transport

Let $M$ be a manifold, $\mathscr{F}(M)$ the set of smooth $\mathbb{R}$-valued scalar functions and $\mathfrak{X}(M)$ the set of real vector fields on $M$. For $t \in \mathbb{R}$ and $f \in \mathfrak{X}(M)$ let $\Psi_{t,f}\colon M \to M$ denote the solution operator such that the differential equation $\gamma'(t) = f(\gamma(t))$, $\gamma(0) = p \in M$ has solution $\gamma(t) = \Psi_{t,f}(p)$. For $\phi \in \mathscr{F}(M)$ we define *pullback along the flow* $\Psi^*_{t,f}\colon \mathscr{F}(M) \to \mathscr{F}(M)$ as

$$\Psi^*_{t,f}\phi = \phi \circ \Psi_{t,f}.$$

The directional derivative $f(\phi) \in \mathscr{F}(M)$ is defined as

$$f(\phi) = \left.\frac{d}{dt}\right|_{t=0} \Psi^*_{t,f}\phi.$$

Through this, we identify $\mathfrak{X}(M)$ with the first order derivations of $\mathscr{F}(M)$, and we obtain higher order derivations by iterating, i.e. $f * f$ is the second order derivation $f * f(\phi) := f(f(\phi))$. With $\mathbb{I}\phi = \phi$ being the 0-order identity operator, the set of all higher order differential operators on $\mathscr{F}(M)$ is called the *universal enveloping algebra* $U(\mathfrak{X}(M))$. This is an algebra with an associative product $*$. The pullback satisfies

$$\frac{d}{dt}\Psi^*_{t,f}\phi = \Psi^*_{t,f}f(\phi).$$

By iteration we find that $\left.\frac{d^n}{dt^n}\right|_{t=0}\Psi^*_{t,f}\phi = f(f(\cdots f(\phi))) = f^{*n}(\phi)$ and hence the Taylor expansion of the pullback is

$$\Psi^*_{t,f}\phi = \phi + tf(\phi) + \frac{t^2}{2!}f * f(\phi) + \cdots = \exp^*(tf)(\phi), \tag{2}$$

where we define the exponential as

$$\exp^*(tf) := \sum_{j=0}^{\infty} \frac{t^j}{j!} f^{*j}.$$

This exponential is an element of $U(\mathfrak{X}(M))$, or more correctly, since it is an infinite series, in the completion of this algebra. We can recover the flow $\Psi_{t,f}$ from $\exp^*(tf)$ by letting $\phi$ be the coordinate maps. However, some caution must be exercised, since pullbacks compose contravariantly $\left(\Psi_{t,f} \circ \Psi_{s,g}\right)^* = \Psi^*_{s,g} \circ \Psi^*_{t,f}$, we have that $\exp^*(sg) * \exp^*(tf)$ corresponds to the diffeomorphism $\Psi_{t,f} \circ \Psi_{s,g}$.

Numerical integrators are constructed by sampling a vector field in points near a base point. To understand this process, we need to transport vector fields. Pullback of vector fields is, however, less canonical than of scalar functions. The differential geometric concept of parallel transport of vectors is defined in terms of a *connection*. An affine connection is a $\mathbb{Z}$-bilinear mapping $\rhd\colon \mathfrak{X}(M)\times\mathfrak{X}(M)\rightarrow\mathfrak{X}(M)$ such that

$$
\begin{aligned}
(\phi f)\rhd g &= \phi(f\rhd g)\\
f\rhd(\phi g) &= f(\phi)g+\phi(f\rhd g)
\end{aligned}
$$

for all $f,g\in\mathfrak{X}(M)$ and $\phi\in\mathscr{F}(M)$. Note that the standard notation for a connection in differential geometry is is $\nabla_f g\equiv f\rhd g$. Our notation is chosen to emphasise the operation as a binary product on the set of vector fields. The triangle notation looks nicer when we iterate, such as in (3) below. Furthermore, the triangle notation is also standard in much of the algebraic literature on pre-Lie algebras, as well as in several recent works on post-Lie algebras.

There is an intimate relationship between connections and the concept of parallel transport. For a curve $\gamma(t)\in M$, let $\Gamma(\gamma)_s^t$ denote *parallel transport along* $\gamma(t)$, meaning that

- $\Gamma(\gamma)_s^t\colon TM_{\gamma(s)}\rightarrow TM_{\gamma(t)}$ is a linear isomorphism of the tangent spaces.
- $\Gamma(\gamma)_s^s=\mathrm{Id}$, the identity map.
- $\Gamma(\gamma)_t^u\circ\Gamma(\gamma)_s^t=\Gamma(\gamma)_s^u$.
- $\Gamma$ depends smoothly on $s,t$ and $\gamma$.

From $\Gamma$, let us consider the action of *parallel transport pullback* of vector fields, for $t\in\mathbb{R}$ and $f\in\mathfrak{X}(M)$ we denote $\Psi_{t,f}^*\colon\mathfrak{X}(M)\rightarrow\mathfrak{X}(M)$ the operation

$$\Psi_{t,f}^* g(p):=\Gamma(\gamma)_t^0 g(\gamma(t)),\quad \text{for the curve } \gamma(t)=\Psi_{t,f}(p).$$

Any connection can be obtained from a parallel transport as the rate of change of the parallel transport pullback. For a given $\Gamma$ we can define a corresponding connection as

$$f\rhd g:=\left.\frac{d}{dt}\right|_{t=0}\Psi_{t,f}^* g.$$

Conversely, we can recover $\Gamma$ from $\rhd$ by solving a differential equation. We seek a power series expansion of the parallel transport pullback. Just like the case of scalars, it holds also for pullback of vector fields that

$$\frac{\partial}{\partial t}\Psi_{t,f}^* g=\Psi_{t,f}^* f\rhd g,$$

hence we obtain the following Taylor series of the pullback

$$\Psi_{t,f}^* g=g+tf\rhd g+\frac{t^2}{2}f\rhd(f\rhd g))+\frac{t^3}{3!}f\rhd(f\rhd(f\rhd g)))+\cdots. \tag{3}$$

Recall that in the case of pullback of a scalar function, we used $f(g(\phi))=(f*g)(\phi)$ to express the pull-back in terms of $\exp^*(tf)$. Whether or not we can do similarly for vector fields depends on geometric properties of the connection. We would like to extend $\rhd$ from $\mathfrak{X}(M)$ to $U(\mathfrak{X}(M))$ such that $f\rhd(g\rhd h)=(f*g)\rhd h$ and hence (3) becomes $\Psi^*_{t,f} g = \exp^*(tf) \rhd g$. However, this requires that $f \rhd (g \rhd h) - g \rhd (f \rhd h) = [\![f,g]\!] \rhd h$, where $[\![f,g]\!] := f * g - g * f$ is the Jacobi bracket of vector fields. The *curvature tensor* of the connection $R\colon \mathfrak{X}(M) \wedge \mathfrak{X}(M) \to \mathrm{End}(\mathfrak{X}(M))$ is defined as

$$R(f,g)h := f \rhd (g \rhd h) - g \rhd (f \rhd h) - [\![f,g]\!] \rhd h.$$

Thus, we only expect to find a suitable extension of $\rhd$ to $U(\mathfrak{X}(M))$ if $\rhd$ is *flat*, i.e. when $R = 0$.

In addition to the curvature, the other important tensor related to a connection is the torsion. Given $\rhd$, we define an $\mathscr{F}(M)$-bilinear mapping $\cdot\colon \mathfrak{X}(M) \times \mathfrak{X}(M) \to U(\mathfrak{X}(M))$ as

$$f \cdot g := f * g - f \rhd g. \tag{4}$$

The skew-symmetrisation of this product called the *torsion*

$$T(f,g) := g \cdot f - f \cdot g \in \mathfrak{X}(M),$$

and if $f \cdot g = g \cdot f$ we say that $\rhd$ is *torsion free*.

The standard connection on $\mathbb{R}^n$ is flat and torsion free. In this case the algebra $\{\mathfrak{X}(M), \rhd\}$ forms a *pre-Lie* algebra (defined below). This gives rise to classical B-series. More generally, transport by left or right multiplication on a Lie group yields a flat connection where the product $\cdot$ is associative, but not commutative. The resulting algebra is called *post-Lie* and the series are called *Lie–Butcher series*. A third important example is the Levi–Civita connection on a symmetric space, where $\cdot$ is a Jordan product, $T = 0$ and $R$ is constant, non-zero. This third case is the subject of forthcoming papers, but will not be discussed here.

## 2.2 The flat Cartan connection on a Lie group

Let $G$ be a Lie group with Lie algebra $\mathfrak{g}$. For $V \in \mathfrak{g}$ and $p \in G$ we let $Vp := TR_pV \in T_pG$. There is a 1–1 correspondence between functions $f \in C^\infty(G,\mathfrak{g})$ and vector fields $\xi_f \in \mathfrak{X}(G)$ given as $\xi_f(p) = f(p)p$. Left multiplication with $q \in G$ gives rise to a parallel transport

$$\Gamma_q\colon T_pG \to T_{qp}G\colon Vp \mapsto Vqp.$$

This transport is independent of the path between $p$ and $qp$ and hence gives rise to a flat connection. We express the corresponding parallel transport pullback on the space $C^\infty(G,\mathfrak{g})$ as

$$(\Gamma_q^* f)(p) = f(qp)$$

which yields the flat connection

$$(f \rhd g)(q) = \left.\frac{d}{dt}\right|_{t=0} g(\exp(tf(q))q).$$

The torsion is given as [22]

$$T(f,g)(p) = -[f(p), g(p)]_{\mathfrak{g}}.$$

The two operations $f \rhd g$ and $[f,g] := -T(f,g)$ turn $C^\infty(G,\mathfrak{g})$ into a *post-Lie algebra*, see Definition 3 below. This is the foundation of Lie–Butcher series.

We can alternatively express the connection and torsion on $\mathfrak{X}(G)$ via a basis $\{E_j\}$ for $\mathfrak{g}$. Let $\partial_j \in \mathfrak{X}(G)$ be the right invariant vector field $\partial_j(p) = E_j p$. For $F, G \in \mathfrak{X}(G)$, where $F = f^i\partial_i$, $G = g^j\partial_j$ [1] and $f^i, g^j \in \mathscr{F}(G)$, we have

$$F \rhd G = f^i \partial_i(g^j)\partial_j$$
$$F \cdot G = f^i g^j \partial_i \partial_j$$
$$T(F,G) = f^i g^j (\partial_i\partial_j - \partial_j\partial_i).$$

We return to $\rhd$ defined on $C^\infty(G,\mathfrak{g})$. Let $U(\mathfrak{g})$ be the span of the basis $\{E_{j_1}E_{j_2}\cdots E_{j_k}\}$, where $E_{j_1}E_{j_2}\cdots E_{j_k} \in U(\mathfrak{g})$ corresponds to the right invariant $k$-th order differential operator $\partial_{j_k}\cdots\partial_{j_2}\partial_{j_1} \in U(\mathfrak{X}(G))$. On $U(\mathfrak{g})$ we have two different associative products, the composition of differential operators $f * g$ and the 'concatenation product' $f \cdot g = f * g - f \rhd g$ which is computed as the concatenation of the basis, $f^i E_i \cdot g^j E_j = f^i g^j E_i E_j$. The general relationship between these two products and $\rhd$ extended to $U(\mathfrak{g})$ is given in (28)–(31) below. In particular we have

$$f \rhd (g \rhd h) = (f * g) \rhd h,$$

which yields the exponential form of the parallel transport

$$\Psi^*_{t,f} g = \exp^*(tf) \rhd g,$$

where $\exp^*(tf)$ is giving us the exact flow of $f$.

We can also form the exponential with respect to the other product,

$$\exp^{\cdot}(tf) = I + tf + \frac{t^2}{2} f \cdot f + \frac{t^3}{3!} f \cdot f \cdot f + \cdots.$$

What is the geometric meaning of this? We say that a vector field $g$ is *parallel along* $f$ if the parallel transport pullback of $g$ along the flow of $f$ is constant, and we say that $g$ is *absolutely parallel* if it is constant under *any* parallel transport. Infinitesimally we have that $g$ is parallel along $f$ if $f \rhd g = 0$ and $g$ is absolutely parallel if $f \rhd g = 0$ for all $f$. In $C^\infty(G,\mathfrak{g})$ the absolutely parallel functions are constants $g(p) = V$, which correspond to right invariant vector fields $\xi_g \in \mathfrak{X}(G)$ given as $\xi_g(p) = Vp$. The flow of parallel vector fields are the geodesics of the connection. If $g$ is absolutely parallel, we have $g * g = g \cdot g + g \rhd g = g \cdot g$, and more generally $g^{n*} = g^{n\cdot}$, hence $\exp^*(g) = \exp^{\cdot}(g)$. If $f(p) = g(p)$ at a point $p \in G$, then they define the same tangent at the point. Hence $f^{n\cdot}(p) = g^{n\cdot}(p)$ for all $n$, and we conclude that $\exp^{\cdot}(f)(p) = \exp^{\cdot}(g)(p) = \exp^*(g)(p)$. Thus, the concatenation exponential $\exp^{\cdot}(f)$ of a general vector field $f$ produces the flow which in a given point follows the geodesic tangent to $f$ at the given point.

---

[1] Einstein summation convention.

On a Lie group, we have for two arbitrary vector fields represented by general functions $f, g \in C^{\infty}(G, \mathfrak{g})$ that

$$(\exp^{\cdot}(tf) \rhd g)(p) = g\,(\exp(tf(p))p)\,. \tag{5}$$

## 2.3 Numerical integration

Lie–Butcher series and its cousins are general mathematical tools with applications in numerics, stochastics and renormalisation. The problem of numerical integration on manifolds is a particular application which has been an important source of inspiration. We discuss a simple illustrative example.

*Example 1 (Lie–trapezoidal method).* Consider the classical trapezoidal method. For a differential equation $y'(t) = f(y(t))$, $y(0) = y_0$ on $\mathbb{R}^n$ a step from $t = 0$ to $t = h$ is given as

$$\begin{aligned} K &= \frac{h}{2}\left(f(y_0) + f(y_1)\right) \\ y_1 &= y_0 + K. \end{aligned}$$

Consider a curve $y(t) \in G$ evolving on a Lie group such that $y'(t) = f(y(t))y(t)$, where $f \in C^{\infty}(G, \mathfrak{g})$ and $y(0) = y_0$. In the Lie-trapezoidal integrator a step from $y_0$ to $y_1 \approx y(h)$ is given as

$$\begin{aligned} K &= \frac{h}{2}\left(f(y_0) + f(y_1)\right) \\ y_1 &= \exp_{\mathfrak{g}}(K)y_0, \end{aligned}$$

where $\exp_{\mathfrak{g}} \colon \mathfrak{g} \to G$ is the classical Lie group exponential. We can write the method as a mapping $\Phi_{\text{trap}} \colon \mathfrak{X}(M) \to \mathrm{Diff}(G)$ from vector fields to diffeomorphisms on $G$, given in terms of parallel transport on $\mathfrak{X}(M)$ as

$$K = \frac{1}{2}\left(f + \exp^{\cdot}(K) \rhd f\right) \tag{6}$$

$$\Phi_{\text{trap}}(f) := \exp^{\cdot}(K). \tag{7}$$

To simplify, we have removed the timestep $h$, but this can be recovered by the substitution $f \mapsto hf$. Note that we present this as a process in $U(\mathfrak{X}(M))$, without a reference to a given base point $y_0$. The method computes a diffeomorphism $\Phi_{\text{trap}}(f)$, which can be evaluated on a given base point $y_0$. This absence of an explicit base point facilitates an interpretation of the method as a process in the enveloping algebra of a free post-Lie algebra, an abstract model of $U(\mathfrak{X}(M))$ to be discussed in the sequel.

A basic problem of numerical integration is to understand in what sense a numerical method $\Phi(tf)$ approximates the exact flow $\exp^*(tf)$. The *order* of the approximation is computed by comparing the LB-series expansion of $\Phi(tf)$ and $\exp^*(tf)$, and comparing to which order in $t$ the two series agree.

The *backward error* of the method is defined as a modified vector field $\widetilde{f_h}$ such that the exact flow of $\widetilde{f_h}$ interpolates the numerical solution at integer times[2]. The combinatorial definition of the backward error is

$$\exp^*(\widetilde{f_h}) = \Phi(hf).$$

The backward error is an important tool which yields important structural information of the numerical flow operator $f \mapsto \Phi(hf)$. The backward error analysis is fundamental in the study of geometric properties of numerical integration algorithms [9, 14].

Yet another problem is the numerical technique of *processing* a vector field, i.e. we seek a modified vector field $\widetilde{f_h}$ such that $\Phi(\widetilde{f_h}) = \exp^*(f)$. An important tool in the analysis of this technique is the characterization of a *substitution law*. What happens to the series expansion of $\Phi(hf)$ if $f$ is replaced by a modified vector field $\widetilde{f_h}$ expressed in terms of a series expansion involving $f$?

The purpose of this essay is not to pursue a detailed discussion of numerical analysis of integration schemes. Instead we want to introduce the algebraic structures needed to formalize the structure of the series expansions. In particular we will present recursive formulas for the basic algebraic operations suitable for computer implementations.

We finally remark that numerical integrators are typically defined as *families of mappings*, given in terms of unspecified coefficients. For example the Runge–Kutta family of integrators can be defined in terms of real coefficients $\{a_{i,j}\}_{i,j=1}^{s}$ and $\{b_j\}_{j=1}^{s}$ as

$$\begin{aligned} K_i &= \exp^{\cdot}\Big(\sum_{j=1}^{s} a_{i,j}K_j\Big) \triangleright f, \quad \text{for } i = 1,\cdots,s \\ \Phi_{\mathrm{RK}}(f) &= \exp^{\cdot}\Big(\sum_{j=1}^{s} b_jK_j\Big). \end{aligned}$$

In a computer package for computing with LB-series we want the possibility of computing series expansions of such parametrized families without specifying the coefficients. This is accomplished by defining the algebraic structures not over the concrete field of real numbers $\mathbb{R}$, but instead allowing this to be replaced by an abstract commutativ ring with unit, such as e.g. the ring of all real polynomials in the indeterminates $\{a_{i,j}\}_{i,j=1}^{s}$ and $\{b_j\}_{j=1}^{s}$.

[2] Technical issues about divergence of the backward error vector field is discussed in [1].

# 3 Algebraic structures of Lie–Butcher theory

We give a concise summary of the basic algebraic structures behind Lie–Butcher series.

## 3.1 Algebras

All vector spaces we consider are over a field[3] k of characteristic 0, e.g. $\mathrm{k} \in \{\mathbb{R}, \mathbb{C}\}$.

**Definition 1 (Algebra).** An algebra $\{\mathscr{A}, *\}$ is a vector space $\mathscr{A}$ with a k-bilinear operation $*\colon \mathscr{A} \times \mathscr{A} \to \mathscr{A}$. $\mathscr{A}$ is called *unital* if it has a unit $\mathbb{I}$ such that $x * \mathbb{I} = \mathbb{I} * x$ for all $x \in \mathscr{A}$. The (minus-)*associator* of the product is defined as

$$a_*(x,y,z) := x*(y*z) - (x*y)*z.$$

If the associator is 0, the algebra is called *associative*.

**Definition 2 (Lie algebra).** A *Lie-algebra* is an algebra $\{\mathfrak{g}, [\cdot,\cdot]\}$ such that

$$[x,y] = -[y,x]$$
$$[[x,y],z] + [[y,z],x] + [[z,x],y] = 0.$$

The bracket $[\cdot,\cdot]$ is called the *commutator* or *Lie bracket*. An associative algebra $\{\mathscr{A}, *\}$ give rise to a Lie algebra $\mathrm{Lie}(\mathscr{A})$, where $[x,y] = x*y - y*x$.

A connection on a fibre bundle which is flat and with constant torsion satisfies the algebraic conditions of a *post-Lie algebra* [22]. This algebraic structure first appeared in a purely operadic setting in [27].

**Definition 3 (Post-Lie algebra).** A *post-Lie algebra* $\{\mathscr{P}, [\cdot,\cdot], \rhd\}$ is a Lie algebra $\{\mathscr{P}, [\cdot,\cdot]\}$ together with a bilinear operation $\rhd\colon \mathscr{P} \times \mathscr{P} \to \mathscr{P}$ such that

$$x \rhd [y,z] = [x \rhd y, z] + [x, y \rhd z] \tag{8}$$
$$[x,y] \rhd z = a_{\rhd}(x,y,z) - a_{\rhd}(y,x,z). \tag{9}$$

A post-Lie algebra defines a relationship between *two* Lie algebras [22].

**Lemma 1.** *For a post-Lie algebra $\mathscr{P}$ the bi-linear operation*

$$[\![x,y]\!] = x \rhd y - y \rhd x + [x,y] \tag{10}$$

*defines another Lie bracket.*

Thus, we have two Lie algebras $\mathfrak{g} = \{\mathscr{P}, [\cdot,\cdot]\}$ and $\overline{\mathfrak{g}} = \{\mathscr{P}, [\![\cdot,\cdot]\!]\}$ related by $\rhd$.

---

[3] In the computer implementations we are relaxing this to allow k more generally to be a commutative ring, such as e.g. polynomials in a set of indeterminates. In this latter case the k-vector space should instead be called a free k-module. We will not pursue this detail in this exposition.

**Definition 4 (Pre-Lie algebra).** A *pre-Lie algebra* $\{\mathscr{L}, \triangleright\}$ is a post-Lie algebra where $[\cdot,\cdot] \equiv 0$, in other words an algebra such that

$$a_{\triangleright}(x,y,z) = a_{\triangleright}(y,x,z).$$

Pre- and post-Lie algebras appear naturally in differential geometry where post-Lie algebras are intimately linked with the differential geometry of Lie groups and pre-Lie algebras with Abelian Lie groups (Euclidean spaces).

## *3.2 Morphisms and free objects*

All algebras of a given type form a *category*, which can be thought of as a directed graph where each node (object) represents an algebra of the given type and the arrows (edges) represent morphisms. Any composition of morphisms is again a morphism. Morphisms are mappings preserving the given algebraic structure. E.g. an algebra morphism $\phi\colon \mathscr{A} \to \mathscr{A}'$ is a k-linear map satisfying $\phi(x * y) = \phi(x) * \phi(y)$. A post-Lie morphism is, similarly, a linear mapping $\phi\colon \mathscr{P} \to \mathscr{P}'$ satisfying both $\phi([x,y]) = [\phi(x), \phi(y)]$ and $\phi(x \triangleright y) = \phi(x) \triangleright \phi(y)$.

In a given category a *free object over a set C* can informally be thought of as a generic algebraic structure. The only equations that hold between elements of the free object are those that follow from the defining axioms of the algebraic structure. Furthermore the free object is not larger than strictly necessary to be generic. Each of the elements of $C$ correspond to generators of the free object. In software a free object can be thought of as a symbolic computing engine; formulas, identities and algebraic simplifications derived within the free object can be applied to any other object in the category. Thus, a detailed understanding of the free objects is crucial for the computer implementation of a given algebraic structure.

**Definition 5 (Free object over a set $C$).** In a given category we define[4] the free object over a set $C$ as an object $\mathrm{Free}(C)$ together with a map $\mathrm{inj}\colon C \hookrightarrow \mathrm{Free}(C)$, called the canonical injection, such that for any object $B$ in the category and any mapping $\phi\colon C \to B$ there exists a unique morphism $!\colon \mathrm{Free}(C) \to B$ such that the diagram commutes

$$\begin{array}{ccc} C & \xrightarrow{\mathrm{inj}} & \mathrm{Free}(C) \\ & \searrow^{\phi} & \downarrow ! \\ & & B \end{array} \quad . \tag{11}$$

We will often consider $C \subset \mathrm{Free}(C)$ without mentioning the map inj.

[4] This definition is not strictly categorical, since the mappings inj and $\phi$ are not morphisms inside a category, but mappings from a set to an object of another category. A proper categorical definition of a free object, found in any book on category theory, is based on a *forgetful functor* mapping the given category into the category of sets. The *free functor* is the left adjoint of the forgetful functor.

*Note 1.* In category theory a free functor is intimately related to a *monad*, a concept which is central in the programming language Haskell. In Haskell the function "inj" is called "return" and the application of ! on $x \in \mathrm{Free}(C)$ is written $x >== \phi$.

A free object can be implemented in different ways, but different implementations are always algebraically isomorphic.

*Example 2.* **Free** $\mathrm{k}$**-vector space** $\mathrm{k}^{(C)}$**:** Consider $C = \{1,2,3,\ldots\}$ and let $\mathrm{inj}(j) = \mathbf{e_j}$ represent a basis for $\mathrm{k}^{(C)}$. Then $\mathrm{k}^{(C)}$ consists of all *finite* $\mathbb{R}$-linear combinations of the basis vectors. Equivalently, we can consider $\mathrm{k}^{(C)}$ as the set of all functions $C \to \mathrm{k}$ with finite support. The unique morphism property states that a linear map is uniquely specified from its values on a set of basis vectors in its domain.

*Example 3.* **Free (associative and unital) algebra** $\mathrm{k}\langle C\rangle$**:** Think of $C$ as an alphabet (collection of letters) $C = \{a,b,c,\ldots\}$. Let $C^*$ denote all *words* over the alphabet, including the empty word $\mathbb{I}$,

$$C^* = \{\mathbb{I}, a, b, c, \ldots, aa, ab, ac, \ldots ba, bb, bc, \ldots\}.$$

Then $\mathrm{k}\langle C\rangle = \{\mathrm{k}^{(C^*)}, \cdot\}$, is the vector space containing finite linear combinations of empty and non-empty words, equipped with a product $\cdot$ which on words is concatenation. Example $aba \cdot cus = abacus$, $\mathbb{I} \cdot abba = abba \cdot \mathbb{I} = abba$. This extends by linearity to $\mathrm{k}^{(C^*)}$ and yields an associative unital algebra. This is also called the non-commutative polynomial ring over $C$.

*Example 4.* **Free Lie algebra** $\mathrm{Lie}(C)$**:** Again, think of $C = \{a,b,c,d,\ldots\}$ as an alphabet. $\mathrm{Lie}(C) \subset \mathrm{k}\langle C\rangle$ is the linear sub space generated by $C$ under the Lie bracket $[w_1, w_2] = w_1 \cdot w_2 - w_2 \cdot w_1$ induced from the product in $\mathrm{k}\langle C\rangle$, thus $c \in C \Rightarrow c \in \mathrm{Lie}(C)$ and $x, y \in \mathrm{Lie}(C) \Rightarrow x \cdot y - y \cdot x \in \mathrm{Lie}(C)$. A basis for $\mathrm{Lie}(C)$ is given by the set of *Lyndon words* [26]. E.g. for $C = \{a,b\}$ the first Lyndon words $a, b, ab, aab, abb$ (up to length 3) represent the commutators

$$\{a, b, [a,b], [a,[a,b]], [[a,b],b], \ldots\}.$$

Computations in a free Lie algebra are important in many applications [21]. Relations such as $[[a,b],c] + [[b,c],a] = [[a,c],b]$ can be computed in $\mathrm{Lie}(C)$ and applied (evaluated) on concrete data in any Lie algebra $\mathfrak{g}$ via the Lie algebra morphism $\mathscr{F}_\phi \colon \mathrm{Lie}(C) \to \mathfrak{g}$, whenever an association of the letters with data in the concrete Lie algebra is provided through a map $\phi \colon C \to \mathfrak{g}$.

*Example 5.* **Free pre-Lie algebra** $\mathrm{preLie}(C)$**:** Consider $C = \{\bullet, \circ, \ldots\}$ as a set of coloured nodes. In many applications $C = \{\bullet\}$, just a single color, and in that case we omit mentioning $C$. A *coloured rooted tree* is a finite connected directed graph where each node (from $C$) has exactly one outgoing edge, except the 'root' node which has no edge out. We illustrate a tree with the root on the bottom and the direction of the edges being down towards the root. Let $T_C$ denote the set of all coloured rooted trees, e.g.

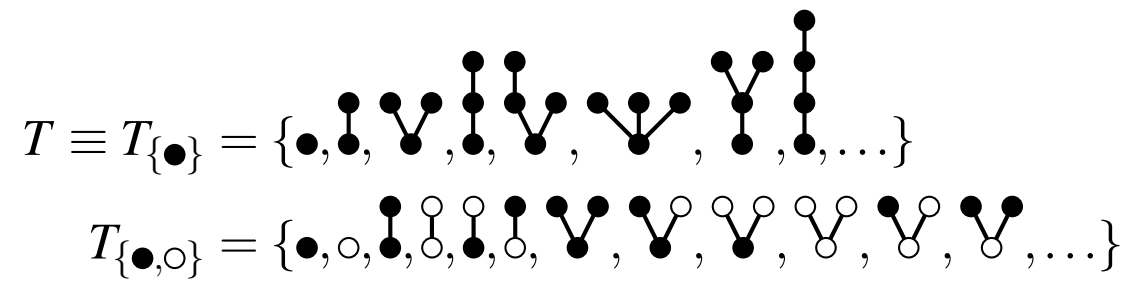

The trees are just graphs without considering an ordering of the branches, so 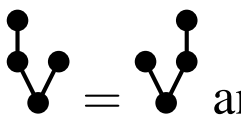 and 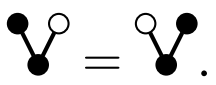. Let $\mathscr{T}_C = \mathrm{k}^{(T_C)}$. The free pre-Lie algebra over $C$ is [6, 10] $\mathrm{preLie}(C) = \{\mathscr{T}_C, \rhd\}$, where $\rhd\colon\ \mathscr{T}_C \times \mathscr{T}_C$ denotes the *grafting product*. For $\tau_1, \tau_2 \in T_C$, the product $\tau_1 \rhd \tau_2$ is the sum of all possible attachments of the root of $\tau_1$ to one of the nodes of $\tau_2$ as shown in this example:

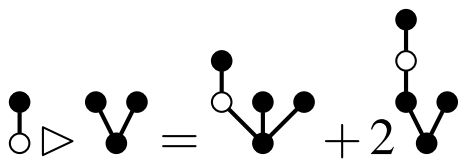

The grafting extends by linearity to all of $\mathscr{T}_C$.

*Example 6.* **Free magma** $\mathrm{Magma}(C) \cong \mathrm{OT}_C$**:** The algebraic definition of a *magma* is a set $C = \{\bullet, \circ, \ldots\}$ with a binary operation $\times$ without any algebraic relations imposed. The free magma over $C$ consists of all possible ways to parenthesize binary operations on $C$, such as $(\bullet \times (\bullet \times \bullet)) \times (\circ \times \bullet)$. There are many isomorphic ways to represent the free magma. For our purpose it is convenient to represent the free magma as ordered (planar[5]) trees with coloured nodes. We let $C$ denote a set of coloured nodes and let $\mathrm{OT}_C$ be the set of all ordered rooted trees with nodes chosen from $C$. On the trees we interpret $\times$ as the *Butcher product* [3]: $\tau_1 \times \tau_2 = \tau$ is a tree where the root of the tree $\tau_2$ is attached to the right part of the root of the tree $\tau_1$, e.g.:

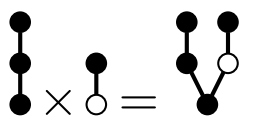

$$= (\bullet \times (\bullet \times \bullet)) \times (\circ \times \bullet).$$

If $C = \{\bullet\}$ has only one element, we write $\mathrm{OT} := \mathrm{OT}_{\{\bullet\}}$. The first few elements of OT are:

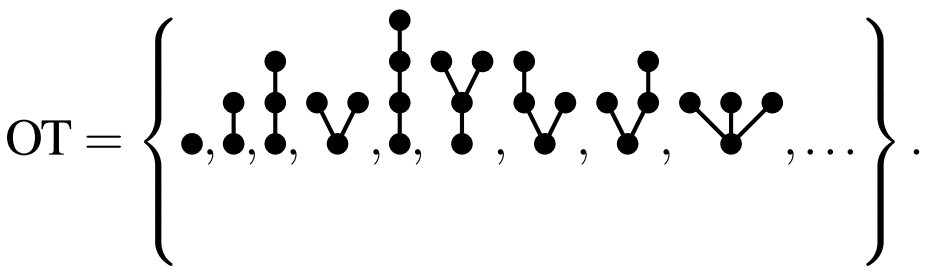

*Example 7.* **The free post-Lie algebra,** $\mathrm{postLie}(C)$, is given as

$$\mathrm{postLie}(C) = \{\mathrm{Lie}(\mathrm{Magma}(C)), \rhd\}, \tag{12}$$

where the product $\rhd$ is defined on $\mathrm{k}^{(\mathrm{Magma}(C))}$ as a derivation of the magmatic product

[5] Trees with different orderings of the branches are considered different, as embedded in the plane.

$$\tau \rhd c = c \times \tau \quad \text{for } c \in C, \tag{13}$$

$$\tau \rhd (\tau_1 \times \tau_2) = (\tau \rhd \tau_1) \times \tau_2 + \tau_1 \times (\tau \rhd \tau_2), \tag{14}$$

and it is extended by linearity and Equations (8)-(9) to all of Lie(Magma($C$)).

Under the identification Magma($C$) $\cong \mathrm{OT}_C$, the product $\rhd\colon\ \mathrm{k}^{(\mathrm{OT}_C)} \times \mathrm{k}^{(\mathrm{OT}_C)} \to \mathrm{k}^{(\mathrm{OT}_C)}$ is given by *left* grafting. For $\tau_1, \tau_2 \in \mathrm{OT}_C$, the product $\tau_1 \rhd \tau_2$ is the sum of all possible attachments of the root of $\tau_1$ to the left side of each node of $\tau_2$ as shown in this example:

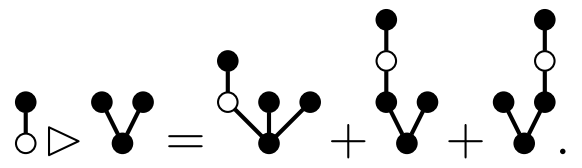

A Lyndon basis for postLie($C$) is given in [20].

## 3.3 Enveloping algebras

Lie algebras, pre- and post-Lie algebras are associated with algebras of first order differential operators (vector fields). Differential operators of higher order are obtained by compositions of these. Algebraically this is described through enveloping algebras.

### 3.3.1 Lie enveloping algebras

Recall that Lie($\cdot$) is a *functor* sending an associative algebra $\mathscr{A}$ to a Lie algebra Lie($\mathscr{A}$), where $[x,y] = x \cdot y - y \cdot x$, and it sends associative algebra homomorphisms to Lie algebra homomorphisms. The universal enveloping algebra of a Lie algebra is defined via a functor $U$ from Lie algebras to associative algebras being the left adjoint of Lie. This means the following:

**Definition 6 (Lie universal enveloping algebra $U(\mathfrak{g})$).** The universal enveloping algebra of a Lie algebra $\mathfrak{g}$ is a unital associative algebra $\{U(\mathfrak{g}), \cdot, \mathbb{I}\}$ together with a Lie algebra morphism inj: $\mathfrak{g} \to \mathrm{Lie}(U(\mathfrak{g}))$ such that for any associative algebra $\mathscr{A}$ and any Lie algebra morphism $\phi\colon \mathfrak{g} \to \mathrm{Lie}(\mathscr{A})$ there exists a unique associative algebra morphism $!\colon U(\mathfrak{g}) \to \mathscr{A}$ such that $\phi = \mathrm{Lie}(!) \circ \mathrm{inj}$.

$$\begin{array}{ccccc} \mathfrak{g} & \xrightarrow{\mathrm{inj}} & \mathrm{Lie}(U(\mathfrak{g})) & \qquad & U(\mathfrak{g}) \\ & \searrow^{\phi} & \downarrow{\scriptstyle \mathrm{Lie}(!)} & & \downarrow{\scriptstyle !} \\ & & \mathrm{Lie}(\mathscr{A}) & & \mathscr{A} \end{array} \tag{15}$$

The *Poincaré–Birkhoff–Witt Theorem* states that for any Lie algebra $\mathfrak{g}$ with a basis $\{e_j\}$, with some total ordering $e_j < e_k$, one gets a basis for $U(\mathfrak{g})$ by taking the set of all *canonical monomials* defined as the non-decreasing products of the basis

elements $\{e_j\}$

$$\text{PBWbasis}(U(\mathfrak{g})) = \{e_{j_1} \cdot e_{j_2} \cdots e_{j_r} \ : \ e_{j_1} \leq e_{j_2} \leq \cdots \leq e_{j_r}, r \in \mathbb{N}\},$$

where we have identified $\mathfrak{g} \subset U(\mathfrak{g})$ using inj. From this it follows that $U(\mathfrak{g})$ is a *filtered* algebra, splitting in a direct sum

$$U(\mathfrak{g}) = \oplus_{j=0}^{\infty} U_j(\mathfrak{g}),$$

where $U_j(\mathfrak{g})$ is the span of the canonical monomials of length $j$, $U_0 = \text{span}(\mathbb{I})$ and $U_1(\mathfrak{g}) \cong \mathfrak{g}$. Furthermore, $U(\mathfrak{g})$ is *connected*, meaning that $U_0 \cong k$, and it is generated by $U_1$, meaning that $U(\mathfrak{g})$ has no proper subalgebra containing $U_1$.

### 3.3.2 Hopf algebras

Recall that a bi-algebra is a unital associative algebra $\{B, \cdot, \mathbb{I}\}$ together with a co-associative co-algebra structure[6] $\{H, \Delta, \varepsilon\}$, where $\Delta\colon B \to B \otimes B$ is the coproduct and $\varepsilon\colon B \to \mathrm{k}$ is the co-unit. The product and coproduct must satisfy the compatibility condition

$$\Delta(x \cdot y) = \Delta(x) \cdot \Delta(y), \tag{16}$$

where the product on the right is componentwise in the tensor product.

**Definition 7 (Hopf algebra).** A *Hopf algebra* $\{H, \cdot, \mathbb{I}, \Delta, \varepsilon, S\}$ is a bi-algebra with an *antipode* $S\colon H \to H$ such that the diagram below commutes.

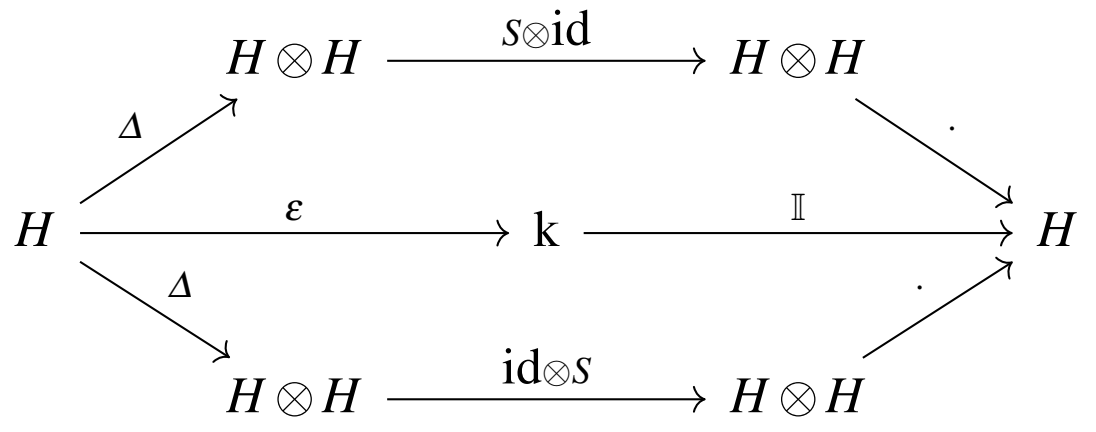

(17)

*Example 8.* **The concatenation de-shuffle Hopf algebra $U(\mathfrak{g})$:** The enveloping algebra $U(\mathfrak{g})$ has the structure of a Hopf algebra, where the coproduct $\Delta_{\sqcup\!\sqcup}\colon U(\mathfrak{g}) \to U(\mathfrak{g}) \otimes U(\mathfrak{g})$ is defined as

$$\Delta_{\sqcup\!\sqcup}(\mathbb{I}) = \mathbb{I} \otimes \mathbb{I} \tag{18}$$

$$\Delta_{\sqcup\!\sqcup}(x) = \mathbb{I} \otimes x + x \otimes \mathbb{I}, \quad \text{for all } x \in \mathfrak{g} \tag{19}$$

$$\Delta_{\sqcup\!\sqcup}(x \cdot y) = \Delta_{\sqcup\!\sqcup}(x) \cdot \Delta_{\sqcup\!\sqcup}(y), \quad \text{for all } x, y \in U(\mathfrak{g}). \tag{20}$$

We call this the de-shuffle coproduct, since it is the dual of the shuffle product. The co-unit is defined as

[6] An associative algebra can be defined by commutative diagrams. The co-algebra structure is obtained by reversing all arrows.

$$\varepsilon(\mathbb{I}) = 1 \tag{21}$$

$$\varepsilon(x) = 0, \quad x \in U_j(\mathfrak{g}), \quad j > 0, \tag{22}$$

and the antipode $S\colon U(\mathfrak{g}) \to U(\mathfrak{g})$ as

$$S(x_1 \cdot x_2 \cdots x_j) = (-1)^j x_j \cdots x_2 \cdot x_1 \quad \text{for all } x_1, \ldots, x_j \in \mathfrak{g}. \tag{23}$$

This turns $U(\mathfrak{g})$ into a filtered, connected, co-commutative Hopf algebra. Connected means that $U_0 \cong \mathrm{k}$ and co-commutative that $\Delta_{\sqcup\!\sqcup}$ satisfies the diagrams of a commutative product, with the arrows reversed. The dual of a commutative product is co-commutative.

The *primitive elements* of a Hopf algebra $H$, defined as

$$\mathrm{Prim}(H) = \{x \in H\colon \Delta(x) = x \otimes \mathbb{I} + \mathbb{I} \otimes x\}$$

form a Lie algebra with $[x,y] = x \cdot y - y \cdot x$. The *Cartier–Milnor–Moore theorem* (CMM) states that if H is a connected, filtered, co-commutative Hopf algebra, then $U(\mathrm{Prim}(H))$ is isomorphic to $H$ as a Hopf algebra. A consequence of CMM is that the enveloping algebra of a free Lie algebra over a set $C$ is given as

$$U(\mathrm{Lie}(C)) = \mathrm{k}\langle C\rangle, \tag{24}$$

the non-commutative polynomials in $C$. Thus, a basis for $U(\mathrm{Lie}(C))$ is given by non-commutative monomials (the empty and non-empty words in $C^*$).

### 3.3.3 Post-Lie enveloping algebras

Enveloping algebras of pre- and post-Lie algebras are discussed by several authors [23, 13, 24, 22]. In our opinion the algebraic structure of the enveloping algebras are easiest to motivate by discussing the post-Lie case, and obtaining the pre-Lie enveloping algebra as a special case. For Lie algebras the enveloping algebras are associative algebras. The corresponding algebraic structure of a post-Lie enveloping algebra is called a D-algebra (D for derivation) [23, 22]:

**Definition 8 (D-algebra).** Let $A$ be a unital associative algebra with a bilinear operation $\rhd\colon A \otimes A \to A$. Write $\mathrm{Der}(A)$ for the set of all $u \in A$ such that $v \mapsto u \rhd v$ is a derivation: $\mathrm{Der}(A) = \{u \in A\colon u \rhd (vw) = (u \rhd v)w + v(u \rhd w) \text{ for all } v, w \in A\}$. We call $A$ a *D-algebra* if for any $u \in \mathrm{Der}(A)$ and any $v, w \in A$ we have

$$\mathbb{I} \rhd v = v \tag{25}$$

$$v \rhd u \in \mathrm{Der}(A) \tag{26}$$

$$(uv) \rhd w = a_{\rhd}(u,v,w) \equiv u \rhd (v \rhd w) - (u \rhd v) \rhd w. \tag{27}$$

In [22] it is shown:

**Proposition 1.** *For any D-algebra A the set of derivations forms a post-Lie algebra*

$$\mathrm{postLie}(A) := \{\mathrm{Der}(A), [\cdot,\cdot], \rhd\},$$

*where* $[x,y] = xy - yx$.

Thus, $\mathrm{postLie}(\cdot)$ is a functor from the category of D-algebras to the category of post-Lie algebras. There is a functor $U(\cdot)$ from post-Lie algebras to D-algebras, which is the left adjoint of $\mathrm{postLie}(\cdot)$. We can define post-Lie enveloping algebras similarly to Definition 6. A direct construction of the post-Lie enveloping algebra is obtained by extending $\rhd$ to the Lie enveloping algebra of the post-Lie algebra [22]:

**Definition 9 (Post-Lie enveloping algebra $U(\mathscr{P})$).** Let $\{\mathscr{P}, [\cdot,\cdot], \rhd\}$ be post-Lie, let $\{U_L, \cdot\} = U(\{\mathscr{P}, [\cdot,\cdot]\})$ be the Lie enveloping algebra and identify $\mathscr{P} \subset U_L$. The post-Lie enveloping algebra $U(\mathscr{P}) = \{U_L, \cdot, \rhd\}$ is defined by extending $\rhd$ from $\mathscr{P}$ to $U_L$ according to

$$\mathbb{I} \rhd v = v \tag{28}$$
$$v \rhd \mathbb{I} = 0 \tag{29}$$
$$u \rhd (vw) = (u \rhd v)w + v(u \rhd w) \tag{30}$$
$$(uv) \rhd w = a_{\rhd}(u,v,w) := u \rhd (v \rhd w) - (u \rhd v) \rhd w \tag{31}$$

for all $u \in \mathscr{P}$ and $v, w \in U_L$. This construction yields $U(\cdot)\colon \mathrm{postLie} \to \text{D-algebra}$ as a left adjoint functor of $\mathrm{postLie}(\cdot)$.

A more detailed understanding of $U(\mathscr{P})$ is obtained by considering its Hopf algebra structures. A Lie enveloping algebra is naturally also a Hopf algebra with the de-shuffle coproduct $\Delta_{\sqcup\sqcup}$. With this coproduct $U(\mathscr{P})$ becomes a graded, connected, co-commutative Hopf algebra where $\mathrm{Der}(U(\mathscr{P})) = \mathrm{Prim}(U(\mathscr{P})) = \mathscr{P}$. Furthermore, the coproduct is compatible with $\rhd$ in the following sense [13]:

$$A \rhd \mathbb{I} = \varepsilon(A)$$
$$\varepsilon(A \rhd B) = \varepsilon(A)\varepsilon(B)$$
$$\Delta_{\sqcup\sqcup}(A \rhd B) = \sum_{\Delta_{\sqcup\sqcup}(A), \Delta_{\sqcup\sqcup}(B)} (A_{(1)} \rhd B_{(1)}) \otimes (A_{(2)} \rhd B_{(2)})$$

for all $A, B \in U(\mathscr{P})$. Here and in the sequel we employ Sweedler's notation for coproducts,

$$\Delta(A) =: \sum_{\Delta(A)} A_{(1)} \otimes A_{(2)}.$$

[7] Sometimes we need a repeated use of a coproduct. Let $\Delta\omega = \sum \omega_{(1)} \otimes \omega_{(2)}$. We continue by using $\Delta$ to split either $\omega_{(1)}$ or $\omega_{(2)}$. Since the coproduct is co-associative this yields the same result $\Delta^2\omega = \sum \omega_{(1)} \otimes \omega_{(2)} \otimes \omega_{(3)}$, and $n$ applications are denoted

[7] Splitting with regard to the coproduct $\Delta$.

$$\Delta^n(A) =: \sum_{\Delta^n(A)} A_{(1)} \otimes A_{(2)} \otimes \cdots \otimes A_{(n+1)}.$$

Just as a post-Lie algebra always has two Lie algebras $\mathfrak{g}$ and $\overline{\mathfrak{g}}$, the post-Lie enveloping algebra $U(\mathscr{P})$ has two associative products $x, y \mapsto xy$ from the enveloping algebra $U(\mathfrak{g})$ and $x, y \mapsto x * y$ from $U(\overline{\mathfrak{g}})$. Both of these products define Hopf algebras with the same unit $\mathbb{I}$, co-unit $\varepsilon$ and de-shuffle coproduct $\Delta_{\sqcup\!\sqcup}$, but with different antipodes.

**Proposition 2.** *[13] On $U(\mathscr{P})$ the product*

$$A * B := \sum_{\Delta_{\sqcup\!\sqcup}(A)} A_{(1)}(A_{(2)} \rhd B) \tag{32}$$

*is associative. Furthermore $\{U(\mathscr{P}), *, \Delta_{\sqcup\!\sqcup}\} \cong U(\overline{\mathfrak{g}})$ are isomorphic as Hopf algebras.*

The following result is crucial for handling the non-commutativity and non-associativity of $\rhd$:

**Proposition 3.** *[23, 13] For all $A, B, C \in U(\mathscr{P})$ we have*

$$A \rhd (B \rhd C) = (A * B) \rhd C. \tag{33}$$

**The free enveloping post-Lie algebra.**

Finally we introduce the enveloping algebra of the free post-Lie algebra $U(\mathrm{postLie}(C))$. Due to CMM, we know that it is constructed from the Hopf algebra

$$U(\mathrm{postLie}(C)) = U(Lie(\mathrm{OT}_C)) = k\langle \mathrm{OT}_C \rangle,$$

i.e. finite linear combinations of *words of ordered trees*, henceforth called *(ordered) forests* $\mathrm{OF}_C$. If $C$ contains only one element, we call the forests OF:

$\mathrm{OF} = \{\mathbb{I},$ 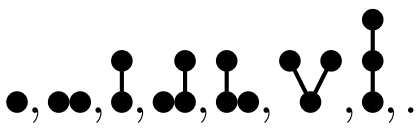 $,\ldots\}$

The Hopf algebra has concatenation of forests as product and coproduct $\Delta_{\sqcup\!\sqcup}$ being de-shuffle of forests. Upon this we define $\rhd$ as left grafting on ordered trees, extended to forests by (28)-(31), where $\mathbb{I}$ is the empty forest, $u$ is an ordered tree and $v, w$ are ordered forests. The left grafting of a forest on another is combinatorially the sum of all possible left attachments of the roots of trees in the left forest to the nodes of the right forest, maintaining order when attaching to the same node, as in this example

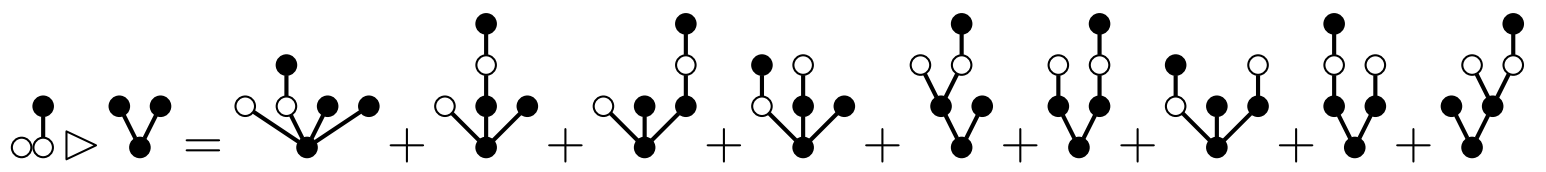

**Four Hopf algebras on ordered forests.**

On $k\langle OF_C\rangle$ we have two associative products $*$ and the concatenation product, denoted $\cdot$. Both these form Hopf algebras with the de-shuffle coproduct $\Delta_{\sqcup\sqcup}$ and antipodes $S_\cdot$ and $S_*$, where

$$S_\cdot(\tau_1\cdot\tau_2\cdots\tau_k) = (-1)^k\tau_k\cdots\tau_2\cdot\tau_1, \quad \text{for } \tau_1\cdot\tau_2\cdots\tau_k\in OT_C$$

and $S_*$ given in (71). With their duals, we have the following four Hopf algebras:

$$\begin{aligned}
\mathscr{H}_\cdot &= \{k\langle OF_C\rangle, \Delta_{\sqcup\sqcup}, \cdot, S_\cdot\}\\
\mathscr{H}_* &= \{k\langle OF_C\rangle, \Delta_{\sqcup\sqcup}, *, S_*\}\\
\mathscr{H}'_\cdot &= \{k\langle OF_C\rangle, \Delta_\cdot, \sqcup\sqcup, S_\cdot\}\\
\mathscr{H}'_* &= \{k\langle OF_C\rangle, \Delta_*, \sqcup\sqcup, S_*\}.
\end{aligned}$$

The four share the same unit $\mathbb{I}\colon\ k\to\mathscr{H}\colon\ 1\mapsto\mathbb{I}$ and the same co-unit $\varepsilon\colon\ \mathscr{H}\to k$, where $\varepsilon(\mathbb{I}) = 1$ and $\varepsilon(\omega) = 0$ for all $\omega\in OF_C\backslash\{\mathbb{I}\}$. All four Hopf algebras are connected and graded with $|\omega|$ counting the number of nodes in a forest. $\mathscr{H}_\cdot$ and $\mathscr{H}'_\cdot$ are also connected and graded with the word length as a grading, although this grading is of less importance for our applications.

### 3.3.4 Lie–Butcher series

The vector space $k\langle OT_C\rangle$ consists of *finite* linear combinations of forests. In order to be able to symbolically represent flow maps and backward error analysis, we do, however, need to extend the space to infinite sums. For a (non-commutative) polynomial ring $k\langle C\rangle$, we denote $k\langle\langle C\rangle\rangle$ the set of infinite (formal) power series. Let $\langle\cdot,\cdot\rangle\colon\ k\langle C\rangle\times k\langle C\rangle\to k$ denote the inner product where the monomials (words in $C^*$) form an orthonormal basis. This extends to a dual pairing

$$\langle\cdot,\cdot\rangle\colon\ k\langle\langle C\rangle\rangle\times k\langle C\rangle\to k, \tag{34}$$

which identifies $k\langle\langle C\rangle\rangle = k\langle C\rangle^*$ as the linar dual space. Any $\alpha\in k\langle\langle C\rangle\rangle$ is uniquely determined by its evaluation on the finite polynomials, and we may write $\alpha$ as a formal infinite sum

$$\alpha = \sum_{w\in C^*}\langle\alpha,w\rangle w.$$

Any k-linear map $f\colon\ k\langle\langle C\rangle\rangle\to k\langle\langle C\rangle\rangle$ can be computed from its dual $f^*\colon\ k\langle C\rangle\to k\langle C\rangle$ as $\langle f(\alpha),w\rangle = \langle\alpha, f^*(w)\rangle$for all $w\in C^*$.

**Definition 10 (Lie–Butcher series $LB(C)$).** The Lie–Butcher series over a set $C$ is defined as

$$LB(C) := U(\text{postLie}(C))^*.$$

This is the vector space $\mathrm{k}\langle\langle \mathrm{OT}_C \rangle\rangle$ (infinite linear combinations of ordered forests). All the operations we consider on this space are defined by their duals acting upon $\mathrm{k}\langle \mathrm{OT}_C \rangle$, see Section 4.1.

The space $\mathrm{LB}(C)$ has two important subsets, the *primitive elements* and the *group like elements*.

**Definition 11 (Primitive elements $\mathfrak{g}_{\mathrm{LB}}$).** The primitive elements of $\mathrm{LB}(C)$, denoted $\mathfrak{g}_{\mathrm{LB}}$ are given as

$$\mathfrak{g}_{\mathrm{LB}} = \{\alpha \in \mathrm{LB}(C) \colon \Delta_{⧢}(\alpha) = \alpha \otimes \mathbb{I} + \mathbb{I} \otimes \alpha\}, \tag{35}$$

where $\Delta_{⧢}$ is the graded completion of the de-shuffle coproduct. This forms a post-Lie algebra which is the graded completion of the free post-Lie algebra $\mathrm{postLie}(C)$.

**Definition 12 (The Lie–Butcher group $G_{\mathrm{LB}}$).** The group like elements of $\mathrm{LB}(C)$, denoted $G_{\mathrm{LB}}$ are given as

$$G_{\mathrm{LB}} = \{\alpha \in \mathrm{LB}(C) \colon \Delta_{⧢}(\alpha) = \alpha \otimes \alpha\}, \tag{36}$$

where $\Delta_{⧢}$ is the graded completion of the de-shuffle coproduct.

The Lie–Butcher group is a group both with respect to the concatenation product and the product $*$ in (32). There are also two exponential maps with respect to the two associative products sending primitive elements to group-like elements

$$\exp, \exp^* \colon \mathfrak{g}_{\mathrm{LB}} \to G_{\mathrm{LB}}.$$

Both these are 1–1 mappings with inverses given by the corresponding logarithms

$$\log, \log^* \colon G_{\mathrm{LB}} \to \mathfrak{g}_{\mathrm{LB}}.$$

# 4 Computing with Lie–Butcher series

In this section, we will list important operations on Lie–Butcher series. A focus will be given on recursive formulations which are suited for computer implementations.

## *4.1 Operations on infinite series computed by dualisation*

Lie–Butcher series are infinite series, and in principle the only computation we consider on an infinite series is the evaluation of the dual pairing (34). All operations on infinite Lie–Butcher series, $\alpha \in \mathrm{LB}(C)$, are computed by dualisation, throwing the operation over to the finite right hand part of the dual pairing. By recursions, the dual computation on the right hand side is moving towards terms with a lower grade, and finally terminates. Some modern programming languages, such as Haskell, allow for *lazy evaluation*, meaning that terms are not computed before they are needed to produce a result. This way it is possible to implement proper infinite series.

*Example 9.* The computation of the de-shuffle coproduct of infinite series can be computed as

$$\langle \Delta_{\sqcup\!\sqcup}(\alpha), \omega_1 \otimes \omega_2 \rangle = \langle \alpha, \omega_1 \sqcup\!\sqcup\, \omega_2 \rangle, \tag{37}$$

where the pairing on the left is defined componentwise in the tensor product,

$$\langle \alpha_1 \otimes \alpha_2, \omega_1 \otimes \omega_2 \rangle = \langle \alpha_1, \omega_1 \rangle \cdot \langle \alpha_2, \omega_2 \rangle$$

and shuffle product $\omega \sqcup\!\sqcup\, \widetilde{\omega}$ of two words in an alphabet is the sum over all permutations of $\omega\widetilde{\omega}$ which are not changing the internal order of the letters coming from each part, e.g.

$$ab \sqcup\!\sqcup\, cd = abcd + acbd + cabd + acdb + cadb + cdab.$$

A recursive formula for the shuffle product is given below.

Any linear operation whose dual sends polynomials in $\mathrm{k}\langle \mathrm{OT}_C \rangle$ to polynomials (or tensor products of these) is well defined on infinite LB-series by such dualisation.

**Linear algebraic operations.**

$$\begin{aligned} +\colon & \ \mathrm{LB}(C) \times \mathrm{LB}(C) \to \mathrm{LB}(C) \quad \text{(addition)} \\ \cdot\colon & \ \mathrm{k} \times \mathrm{LB}(C) \to \mathrm{LB}(C) \quad \text{(scalar multiplication).} \end{aligned}$$

These are computed as $\langle \alpha + \beta, w \rangle = \langle \alpha, w \rangle + \langle \beta, w \rangle$ and $\langle c \cdot \alpha, w \rangle = c \cdot \langle \alpha, w \rangle$. Note that $\mathfrak{g}_{\mathrm{LB}} \subset \mathrm{LB}(C)$ is a linear subspace closed under these operations, $G_{\mathrm{LB}} \subset \mathrm{LB}(C)$ is *not* a linear subspace.

## 4.2 Operations on forests computed by recursions in a magma

Similar to the case of trees, Section 3.2, many recursion formulas for forests are suitably formulated in terms of magmatic products on forests. Let $B^- \colon \mathrm{OT}_C \to \mathrm{OF}_C$ denote the removal of the root, sending a tree to the forest containing the branches of the root, and for every $c \in C$ define $B_c^+ \colon \mathrm{OF}_C \to \mathrm{OT}_C$ as the addition of a root of colour $c$ to a forest, producing a tree, example

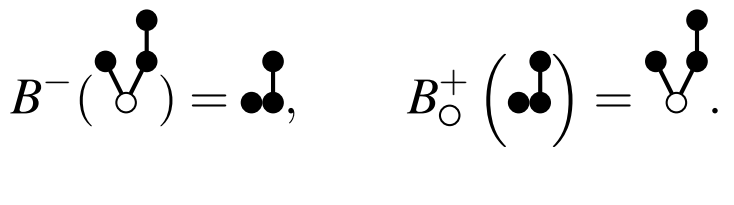

**Definition 13 (Magmatic products on $\mathrm{OF}_C$).** For every $c \in C$, define a product $\times_c \colon \mathrm{OF}_C \times \mathrm{OF}_C \to \mathrm{OF}_C$ as

$$\omega_1 \times_c \omega_2 := \omega_1 B_c^+(\omega_2). \tag{38}$$

In the special case where $C = \{\bullet\}$ contains just one element, then $B^+ \colon \mathrm{OF} \to \mathrm{OT}$ is 1–1, sending the above product on forests to the Butcher product on trees; $B^+(\omega_1 \times_\bullet \omega_2) = B^+(\omega_1) \times B^+(\omega_2)$. Thus, in this case $\{\mathrm{OF}, \times_\bullet\} \cong \{\mathrm{OT}, \times\} \cong \mathrm{Magma}(\{\bullet\})$.

For a general $C$ we have that any $\omega \in \mathrm{OF}_C \backslash \mathbb{I}$ has a unique decomposition

$$\omega = \omega_L \times_c \omega_R, \qquad c \in C, \quad \omega_L, \omega_R \in \mathrm{OF}_C. \tag{39}$$

The set of forests $\mathrm{OF}_C$ is freely generated from $\mathbb{I}$ by these products, e.g.

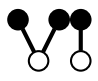 $= (\mathbb{I} \times_\circ ((\mathbb{I} \times_\bullet \mathbb{I}) \times_\bullet \mathbb{I})) \times_\circ (\mathbb{I} \times_\bullet \mathbb{I}).$

Thus, there is a 1–1 correspondence between $\mathrm{OF}_C$ and binary trees where the internal nodes are coloured with $C$. We may take the binary tree representation as the *definition* of $\mathrm{OF}_C$ and express any computation in terms of this.

**Definition 14 (Magmatic definition of $\mathrm{OF}_C$).** Given a set $C$, the ordered forests $\mathrm{OF}_C$ are defined recursively as

$$\mathbb{I} \in \mathrm{OF}_C \tag{40}$$

$$\omega = \omega_L \times_C \omega_R \in \mathrm{OF}_C \quad \text{for evey } \omega_L, \omega_R \in \mathrm{OF}_C \text{ and } c \in C. \tag{41}$$

$\mathrm{OF}_C$ has the following operations:

isEmpty: $\mathrm{OF}_C \to$ bool, defined by isEmpty$(\mathbb{I}) =$ 'true', otherwise 'false'.
Left: $\mathrm{OF}_C \to \mathrm{OF}_C$, defined by $\mathrm{Left}(\omega_L \times_c \omega_R) = \omega_L$.
Right: $\mathrm{OF}_C \to \mathrm{OF}_C$, defined by $\mathrm{Right}(\omega_L \times_c \omega_R) = \omega_R$.
Root: $\mathrm{OF}_C \to C$, defined by $\mathrm{Root}(\omega_L \times_c \omega_R) = c$.

$\mathrm{Left}(\mathbb{I})$, $\mathrm{Right}(\mathbb{I})$ and $\mathrm{Root}(\mathbb{I})$ are undefined.

Any operation on forests can be expressed in terms of these. We can define ordered trees as the subset $\mathrm{OT}_C \subset \mathrm{OF}_C$

$$\mathrm{OT}_C := \{\tau \in \mathrm{OF}_C \colon \mathrm{Left}(\tau) = \mathbb{I}\},$$

and in particular the nodes $C \subset \mathrm{OF}_C$ are identified as $C \cong \{\mathbb{I} \times_c \mathbb{I}\}$. From this we define $B^- \colon \mathrm{OT}_C \to \mathrm{OF}_C$ and $B_c^+ \colon \mathrm{OF}_C \to \mathrm{OT}_C$ as

$$B^-(\tau) = \mathrm{Right}(\tau) \tag{42}$$

$$B_c^+(\omega) = \mathbb{I} \times_c \omega. \tag{43}$$

The Butcher product of two trees $\tau, \tau' \in \mathrm{OT}_C$, where $c = \mathrm{Root}(\tau)$, $c' = \mathrm{Root}(\tau')$ is

$$\tau \times \tau' := B_c^+(B^-(\tau) \times_{c'} B^-(\tau')).$$

## 4.3 Combinatorial functions on ordered forests.

The *order* of $\omega \in \mathrm{OF}_C$, denoted $|\omega| \in \mathbb{N}$, counts the number of nodes in the forest. It is computed by the recursion

$$|\mathbb{I}| = 0 \tag{44}$$

$$|\omega_L \times_\bullet \omega_R| = |\omega_L| + |\omega_R| + 1. \tag{45}$$

This counts the number of nodes in $\omega$.

The *ordered forest factorial*, denoted $\omega\text{¡} \in \mathbb{N}$ is defined by the recursion

$$\mathbb{I}\text{¡} = 1 \tag{46}$$

$$\omega\text{¡} = (\omega_L \times_\bullet \omega_R)\text{¡} = |\omega| \cdot \omega_L\text{¡} \cdot \omega_R\text{¡}. \tag{47}$$

We will see that the ordered factorial is important for characterising the flow map (exact solution) of a differential equation. This is a generalisation of the more well-known *tree factorial function for un-ordered trees*, which is denoted $\tau!$ and defined by the recursion

$$\bullet! = 1$$
$$\tau! = |\tau| \cdot \tau_1! \cdot \tau_2! \cdots \tau_p!$$

for $\tau = B^+(\tau_1\tau_2\cdots\tau_p)$.

The relationship between the classical (unordered) and the ordered tree factorial functions is

$$\sigma(\tau) \sum_{\tau' \sim \tau} \frac{1}{\tau'\text{¡}} = \frac{1}{\tau!},$$

where the sum runs over all ordered trees that are equivalent under permutation of the branches and $\sigma(\tau)$ is the symmetry factor of the tree. This identity can be derived from the relationship between classical B-series and LB-series discussed in Section 4.1 of [23], by comparing the exact flow maps $\exp^*(\bullet)$ in the two cases. We omit details.

*Example 10.*

$1/$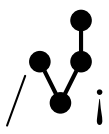¡ $+ 1/$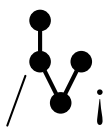¡ $= \frac{1}{12} + \frac{1}{24} = \frac{1}{8} = 1/$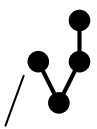$!$

and

$2\Big(1/$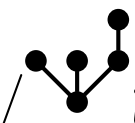¡ $+ 1/$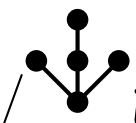¡ $+ 1/$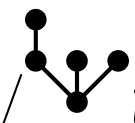¡$\Big) = 2\left(\frac{1}{40} + \frac{1}{60} + \frac{1}{120}\right) = \frac{1}{10} = 1/$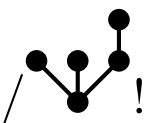$!$.

For the tall tree $\tau = \mathbb{I} \times_{\bullet} (\mathbb{I} \times_{\bullet} (\mathbb{I} \times_{\bullet} (\cdots \times_{\bullet} (\mathbb{I} \times_{\bullet} \mathbb{I}))))$ we have $\tau¡ = \tau! = |\tau|!$. Table 1 on p.35 contains the ordered forest factorial for all ordered forests up to and including order 5.

### 4.4 Concatenation and de-concatenation

Concatenation and de-concatenation

$$\begin{aligned} \cdot :\ & \mathrm{k}\langle \mathrm{OT}_C \rangle \otimes \mathrm{k}\langle \mathrm{OT}_C \rangle \to \mathrm{k}\langle \mathrm{OT}_C \rangle \\ \Delta_{\cdot} :\ & \mathrm{k}\langle \mathrm{OT}_C \rangle \to \mathrm{k}\langle \mathrm{OT}_C \rangle \otimes \mathrm{k}\langle \mathrm{OT}_C \rangle \end{aligned}$$

form a pair of dual operations, just like $\sqcup\!\sqcup$ and $\Delta_{\sqcup\!\sqcup}$ in (37). On monomials $\omega \in \mathrm{OF}_C$ these are given by

$$\begin{aligned} \omega \cdot \omega' &= \omega\omega' \\ \Delta_{\cdot}(\omega) &= \sum_{\substack{\omega_1,\omega_2 \in \mathrm{OF}_C \\ \omega_1 \cdot \omega_2 = \omega}} \omega_1 \otimes \omega_2, \end{aligned}$$

thus for $\omega = \tau_1 \tau_2 \cdots \tau_k$, $\tau_1, \ldots, \tau_k \in \mathrm{OT}_C$ we have

$$\Delta_{\cdot}(\omega) = \omega \otimes \mathbb{I} + \mathbb{I} \otimes \omega + \sum_{j=1}^{k} \tau_1 \cdots \tau_j \otimes \tau_{j+1} \cdots \tau_k.$$

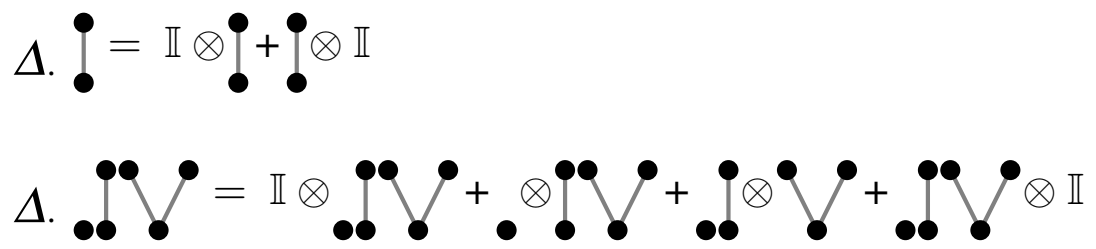

Recursive formulas, where $\widetilde{\omega} \in \mathrm{OF}_C$, $\omega = \omega_L \times_c \omega_R$ are

$$\widetilde{\omega} \cdot \mathbb{I} = \widetilde{\omega} \tag{48}$$

$$\widetilde{\omega} \cdot \omega = (\widetilde{\omega} \cdot \omega_L) \times_c \omega_R \tag{49}$$

and

$$\Delta_{\cdot}(\mathbb{I}) = \mathbb{I} \otimes \mathbb{I} \tag{50}$$

$$\Delta_{\cdot}(\omega) = \Delta_{\cdot}(\omega_L) \cdot (\mathbb{I} \otimes (\mathbb{I} \times_c \omega_R)) + \omega \otimes \mathbb{I}. \tag{51}$$

See Table 2 [8] on p.36 for deconcatenation of all ordered forests up to and including order 4.

The concatenation antipode $S_{\cdot}$, defined in (23), is computed by the recursion

[8] Note that the number under the terms are the coefficients to the terms.

$$S.(\mathbb{I}) = \mathbb{I} \tag{52}$$

$$S.(\omega_L \times_c \omega_R) = -B_c^+(\omega_R) \cdot S.(\omega_L). \tag{53}$$

$S.$ reverse the order of the trees in the forest and negate if there is a odd number of trees in the the forest. See Table 2 on p.36.

### 4.5 *Shuffle and de-shuffle.*

The duality of $\Delta_{\sqcup\!\sqcup}$ and $\sqcup\!\sqcup$ is given in (37). A recursive formula for $\omega \sqcup\!\sqcup \widetilde{\omega}$ where $\omega, \widetilde{\omega} \in \mathrm{OF}_C$ is obtained from the decomposition $\omega = \omega_L \times_c \omega_R$, $\widetilde{\omega} = \widetilde{\omega}_L \times_{\widetilde{c}} \widetilde{\omega}_R$ as

$$\mathbb{I} \sqcup\!\sqcup \omega = \omega \sqcup\!\sqcup \mathbb{I} = \omega \tag{54}$$

$$\omega \sqcup\!\sqcup \widetilde{\omega} = (\omega_L \sqcup\!\sqcup \widetilde{\omega}) \times_c \omega_R + (\omega \sqcup\!\sqcup \widetilde{\omega}_L) \times_{\widetilde{c}} \widetilde{\omega}_R, \tag{55}$$

while (18)-(20) yields the recursion

$$\Delta_{\sqcup\!\sqcup}(\mathbb{I}) = \mathbb{I} \otimes \mathbb{I} \tag{56}$$

$$\Delta_{\sqcup\!\sqcup}(\omega) = \Delta_{\sqcup\!\sqcup}(\omega_L) \cdot ((\mathbb{I} \times_c \omega_R) \otimes \mathbb{I} + \mathbb{I} \otimes (\mathbb{I} \times_c \omega_R))). \tag{57}$$

The shuffle product $\sqcup\!\sqcup$ of two forests is the summation over all permutations of the trees in the forests while preserving the ordering of the trees in each of the initial forests.

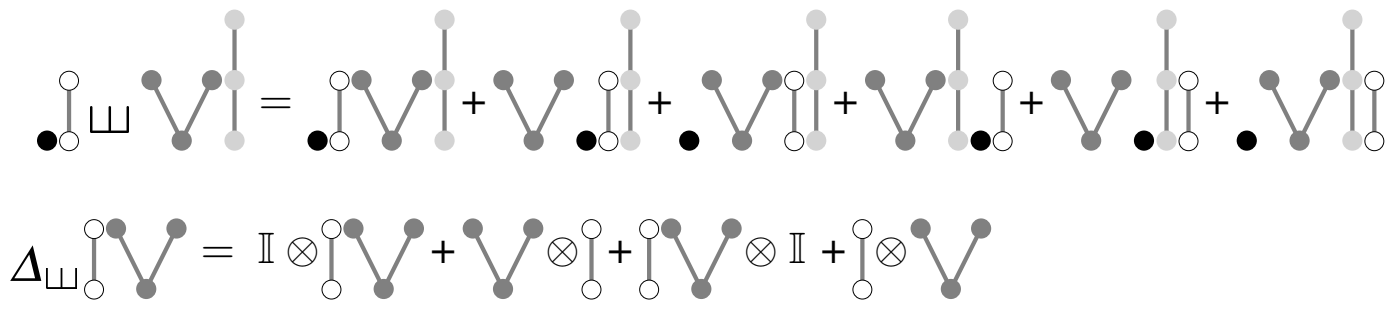

Fig. 1: See Table 2 on p.36 for more examples on deshuffle.

### 4.6 *Grafting, pruning, GL product and GL coproduct*

These are four closely related operations. Grafting is defined in (13)-(14) for trees and (28)-(31) for forests (here $u$ is a tree). Grafting can also be expressed directly through the magmatic definition of $\mathrm{OF}_C$. First we need to decompose $\omega \in \mathrm{OF}_C \backslash \mathbb{I}$ as a concatenation of a tree on the left with a forest on the right, $\omega = \tau' \cdot \omega'$. We define the decomposition $\tau' = \mathrm{LeftTree}(\omega)$, $\omega' = \mathrm{RightForest}(\omega)$ through the following recursions, where $\tau \in \mathrm{OT}_C$ and $\omega = \omega_L \times_c \omega_R$:

$$\mathrm{LeftTree}(\tau) = \tau \tag{58}$$
$$\mathrm{LeftTree}(\omega) = \mathrm{LeftTree}(\omega_L) \tag{59}$$
$$\mathrm{RightForest}(\tau) = \mathbb{I} \tag{60}$$
$$\mathrm{RightForest}(\omega) = \mathrm{RightForest}(\omega_L) \times_c \omega_R. \tag{61}$$

The general recursion for grafting of forests becomes

$$\mathbb{I} \rhd \omega = \omega \tag{62}$$
$$\tau \rhd \mathbb{I} = 0 \tag{63}$$
$$\tau \rhd (\omega_L \times_c \omega_R) = (\tau \rhd \omega_L) \times_c \omega_R + \omega_L \times_c (\tau \cdot \omega_R + \tau \rhd \omega_R) \tag{64}$$
$$(\tau \cdot \omega) \rhd \widetilde{\omega} = \tau \rhd (\omega \rhd \widetilde{\omega}) - (\tau \rhd \omega) \rhd \widetilde{\omega}, \tag{65}$$

for all $\tau \in \mathrm{OT}_C$, $\omega, \widetilde{\omega}, \omega_L, \omega_R \in \mathrm{OF}_C$, $c \in C$. See Table 3 on p.37 for examples.

The associative product $*$ defined in (32) is, in the context of polynomials of ordered trees $\mathrm{k}\langle \mathrm{OT}_C \rangle$, called the (ordered) Grossman–Larsson product [23], *GL product* for short. On $\mathrm{k}\langle \mathrm{OT}_C \rangle$ (and even on $\mathrm{LB}(C)$), we can compute $*$ from grafting as

$$\omega_1 * \omega_2 = B^-(\omega_1 \rhd B^+(\omega_2)).$$

The colour of the added root is irrelevant, since this root is later removed by $B^-$. See Table 3 on p.37 for examples.

The dual of $*$, the GL coproduct $\Delta_* \colon \mathrm{k}\langle \mathrm{OF}_C \rangle \to \mathrm{k}\langle \mathrm{OF}_C \rangle \otimes \mathrm{k}\langle \mathrm{OF}_C \rangle$ has several different characterisations, in terms of left admissible cuts of trees and by recursion [23]. For $\omega = \omega_L \times_c \omega_R$ the recursion is

$$\Delta_*(\mathbb{I}) = \mathbb{I} \otimes \mathbb{I} \tag{66}$$
$$\Delta_*(\omega) = \omega \otimes \mathbb{I} + \Delta_*(\omega_L) \sqcup \times_c \Delta_*(\omega_R), \tag{67}$$

where $\sqcup\times_c \colon \mathrm{k}\langle \mathrm{OT}_C \rangle \otimes \mathrm{k}\langle \mathrm{OT}_C \rangle \otimes \mathrm{k}\langle \mathrm{OT}_C \rangle \otimes \mathrm{k}\langle \mathrm{OT}_C \rangle \to \mathrm{k}\langle \mathrm{OT}_C \rangle \otimes \mathrm{k}\langle \mathrm{OT}_C \rangle$ denotes

$$(\alpha \otimes \widetilde{\alpha}) \sqcup \times_c (\omega \otimes \widetilde{\omega}) := (\alpha \sqcup \omega) \otimes (\widetilde{\alpha} \times_c \widetilde{\omega}).$$

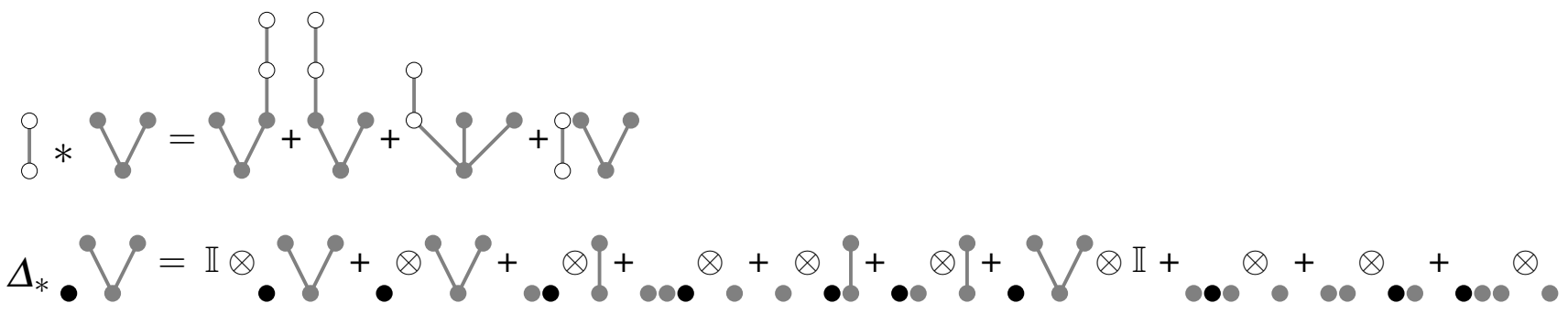

Fig. 2: See Table 3 on p.37 and Table 4 on p.38 for more examples.

The grafting operation $\rhd \colon \mathrm{k}\langle \mathrm{OT}_C \rangle \times \mathrm{k}\langle \mathrm{OT}_C \rangle \to \mathrm{k}\langle \mathrm{OT}_C \rangle$ has a right sided dual we call *pruning*, $\Delta_\rhd \colon \mathrm{k}\langle \mathrm{OT}_C \rangle \to \mathrm{k}\langle \mathrm{OT}_C \rangle \times \mathrm{k}\langle \mathrm{OT}_C \rangle$, dual in the usual sense

$$\langle \alpha \rhd \beta, \omega \rangle = \langle \alpha \otimes \beta, \Delta_{\rhd}(\omega) \rangle.$$

The pruning is characterised by admissible cuts in [17], or it can be computed by the following recursion involving both itself and the GL coproduct,

$$\Delta_{\rhd}(\mathbb{I}) = \mathbb{I} \otimes \mathbb{I} \tag{68}$$

$$\Delta_{\rhd}(\omega_L \times_c \omega_R) = \Delta_{\rhd}(\omega_L) \sqcup\!\sqcup \times_c \Delta_*(\omega_R). \tag{69}$$

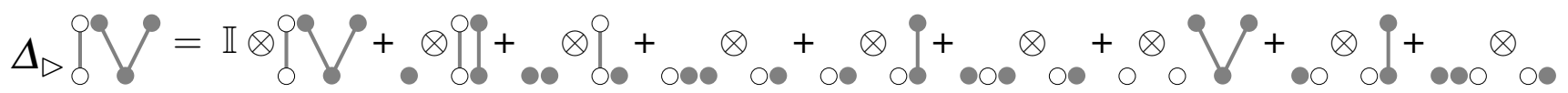

Fig. 3: See also Table 4 on p.38.

**The Lie–Butcher group and the antipode $S_*$.**

The product in the Lie-Butcher group $G_{\mathrm{LB}}$ is the GL product $\alpha, \beta \mapsto \alpha * \beta$. The inverse is given by the *antipode* (with respect to $*$-product), an endomorphism $S_* \in \mathrm{End}(\mathrm{k}\langle \mathrm{OT}_C \rangle)$ such that

$$\langle \alpha^{*-1}, \omega \rangle = \langle \alpha, S_*(\omega) \rangle. \tag{70}$$

A recursive formula for $S_*$ is found in [23]. In our magmatic representation of forests we have

$$S_*(\omega_L \times_c \omega_R) = -\sqcup\!\sqcup\left((S_* \otimes I)(\Delta_*(\omega_L) \sqcup\!\sqcup \times_c \Delta_*(\omega_R))\right). \tag{71}$$

Table 5 on p.39 contain the the result of applying $S_*$ to all ordered forests up to and including order 4.

## 4.7 Substitution, co-substitution, scaling and derivation.

A LB-series is an infinite series of forests built from nodes. The substitution law [7, 8, 17, 4] expresses the operation of replacing each node with an entire LB series. Since a node represents a primitive element, it is necessary to require that the LB-series in the substitution must be an element of $\mathfrak{g}_{\mathrm{LB}}$. The universal property of the free enveloping algebra $U(\mathrm{postLie}(C))$ implies that for any mapping $a\colon C \to \mathscr{P}$ from $C$ into a post-Lie algebra $\mathscr{P}$, there exists a unique D-algebra morphism $!\colon U(\mathrm{postLie}(C)) \to U(\mathscr{P})$ such that the diagram commutes

$$\begin{array}{ccc} C & \xrightarrow{\text{inj}} & U(\mathrm{postLie}(C)) \\ \Big\downarrow{\scriptstyle a} & & \Big\downarrow{\scriptstyle !} \\ \mathscr{P} & \xrightarrow{\quad\text{inj}\quad} & U(\mathscr{P}) \end{array} \tag{72}$$

In particular this holds if $\mathscr{P} = \mathrm{postLie}(C)$, and it also holds if $U(\mathrm{postLie}(C))$ is replaced with its graded completion $\mathrm{LB}(C)$. From this we obtain the algebraic definition of substitution:

**Definition 15 (Substitution).** Given a mapping $a\colon C \to \mathfrak{g}_{\mathrm{LB}}$ there exists a unique D-algebra automorphism $a\star\colon \mathrm{LB}(C) \to \mathrm{LB}(C)$ such that the diagram commutes

$$\begin{array}{ccc} C & \xrightarrow{\text{inj}} & \mathrm{LB}(C) \\ \Big\downarrow a & & \Big\downarrow a\star \\ \mathfrak{g}_{\mathrm{LB}} & \xrightarrow{\text{inj}} & \mathrm{LB}(C). \end{array} \tag{73}$$

This morphism is called *substitution*.

The automorphism property implies that it enjoys many identities such as

$$a\star\mathbb{I} = \mathbb{I} \tag{74}$$
$$a\star(\omega\omega') = (a\star\omega)(a\star\omega') \tag{75}$$
$$a\star(\omega \rhd \omega') = (a\star\omega) \rhd (a\star\omega') \tag{76}$$
$$a\star(\omega * \omega') = (a\star\omega) * (a\star\omega') \tag{77}$$
$$(a\star\otimes a\star)(\Delta_{\sqcup\!\sqcup}(\omega)) = \Delta_{\sqcup\!\sqcup}(a\star\omega). \tag{78}$$

For more details, see [17].

As explained earlier, computations with LB-series are done by considering the series together with a pairing on the space of finite series and computations are performed by deriving how the given operation is expressed as an operation on finite series, via the dual. Thus, to compute substitution of infinite series, we need to characterise the dual map, called co-substitution.

**Definition 16 (Co-substitution).** Given a substitution $a\star\colon \mathrm{LB}(C) \to \mathrm{LB}(C)$, the *co-substitution* $a_\star^T$ is a k-linear map $a_\star^T\colon\ \mathrm{k}\langle \mathrm{OT}_C\rangle \to \mathrm{k}\langle \mathrm{OT}_C\rangle$ such that

$$\langle a\star\beta, x\rangle = \langle \beta, a_\star^T(x)\rangle$$

for all $\beta \in \mathrm{LB}(C)$ and $x \in \mathrm{k}\langle \mathrm{OT}_C\rangle$.

A recursive formula for the co-substitution is derived in [17] in the case where $C = \{\bullet\}$. A general formula for arbitrary finite $C$ is given here, the proof of this formula is similar to the proof in [17] but we omit it. The general formula for $a_\star^T(\omega)$ is based on decomposing $\omega$ with the de-concatenation coproduct $\Delta_\cdot$ and thereafter decomposing the second component with the pruning coproduct $\Delta_\rhd$. To clarify the notation, the decomposition is as follows

$$(I \otimes \Delta_\rhd) \circ \Delta_\cdot(\omega) = \sum_{\Delta_\cdot(\omega)} \sum_{\Delta_\rhd(\omega_{(2)})} \omega_{(1)} \otimes \omega_{(2)(1)} \otimes \omega_{(2)(2)}.$$

With this decomposition, a recursion for $a_\star^T$ is given as $a_\star^T(\mathbb{I}) = \mathbb{I}$ and for $\omega \in \mathrm{OF}_C \backslash \mathbb{I}$

$$a_\star^T(\omega) = \sum_{c\in C}\sum_{\Delta_\cdot(\omega)}\sum_{\Delta_\triangleright(\omega_{(2)})} \left(a_\star^T(\omega_{(1)}) \times_c a_\star^T(\omega_{(2)(1)})\right)\langle a(c), \omega_{(2)(2)}\rangle. \tag{79}$$

The recursion is written more compactly as

$$a_\star^T = \sum_{c\in C} \mu_\cdot \circ (\mu_{\times_c} \otimes I) \circ (a_\star^T \otimes a_\star^T \otimes a(c)) \circ (I \otimes \Delta_\triangleright) \circ \Delta_\cdot\ ,$$

where $\mu_\cdot(\omega \otimes \omega') := \omega \cdot \omega'$, $\mu_{\times_c}(\omega \otimes \omega') := \omega \times_c \omega'$ and $a(c)\colon\ \mathrm{k}\langle \mathrm{OT}_C\rangle \to \mathrm{k}$ denotes $\omega \mapsto \langle a(c), \omega\rangle$.

See Table 6 on p.40 where cosubstitution is calculated for all forests up to and including order 4, assuming $a$ is a infinitesimal character.

Since $a\star$ is compatible with $\Delta_\sqcup$ in the sense of (78), it follows that $a_\star^T$ is a shuffle homomorphism (a character) satisfying

$$a_\star^T(\omega \sqcup \omega') = a_\star^T(\omega) \sqcup a_\star^T(\omega').$$

**Definition 17 (Scaling).** For $t \in \mathrm{k}$ define the map $t(c) = tc\colon C \to \mathfrak{g}_{\mathrm{LB}}$. The corresponding substitution $\alpha \mapsto t \star \alpha$ is called *scaling by t*. For a fixed alpha $t \mapsto t \star \alpha$ defines a curve in $\mathrm{LB}(C)$

Note that $t \star \omega = t^{|\omega|}\omega$ and hence $\langle t \star \alpha, \omega\rangle = t^{|\omega|}\langle \alpha, \omega\rangle$ for all $\omega \in \mathrm{OF}_C$.

**Definition 18 (Derivation).** The derivative of a LB-series $\alpha$, denoted $D\alpha$ is defined as

$$\langle D\alpha, \omega\rangle = |\omega|\langle \alpha, \omega\rangle.$$

Note that if $\mathrm{k} = \mathbb{R}$ we have $D\alpha = \frac{d}{dt}\big|_{t=1}(t \star \alpha)$.

## *4.8 Exponentials and logarithms.*

We have three types of exponential type mappings $\exp^\cdot, \exp^*, \mathrm{evol}\colon \mathfrak{g}_{\mathrm{LB}} \to G_{\mathrm{LB}}$. These are all 1–1 mappings with an inverse being a kind of logarithm. In the interpretation of vector fields on Lie groups, $\exp^\cdot$ defines the geodesics of the connection and $\exp^*$ computes the exact flow of a vector field. The third of these, evol, computes a curve in a Lie group from its development in the Lie algebra i.e. solves an equation of Lie type $y'(t) = y(t)\gamma(t)$ where $\gamma(t) = y^{-1}(t)y'(t)$ is the development of $y(t)$ (left logarithmic derivative). We will have a closer look at these three maps and their inverses.

**Definition 19 (Concatenation exponential).** The *concatenation exponential* $\exp^\cdot\colon \mathfrak{g}_{\mathrm{LB}} \to G_{\mathrm{LB}}$ is defined as

$$\exp^\cdot(\alpha) = \mathbb{I} + \alpha + \frac{1}{2}\alpha\alpha + \frac{1}{6}\alpha\alpha\alpha + \cdots = \sum_{j=0}^{\infty} \frac{1}{j!}\alpha^{\cdot j}. \tag{80}$$

In the algebra $U(\mathrm{postLie}(C))$, with the grading given by PBW, $U_0 = \mathrm{k}\mathbb{I}$, $U_1 = \mathrm{postLie}(C)$ and $U_\ell$ is generated from $U_1$ by $\ell$-fold shuffle products. Since $\langle \exp^{\cdot}(\alpha), x \sqcup\!\sqcup y\rangle = \langle \exp^{\cdot}(\alpha), x\rangle\langle \exp^{\cdot}(\alpha), y\rangle$ we have the following result.

**Lemma 2.** *For $\alpha \in \mathfrak{g}_{\mathrm{LB}}$, the concatenation exponential $\exp^{\cdot}(\alpha)$ is the unique element of $G_{\mathrm{LB}}$ such that $\langle \exp^{\cdot}(\alpha), x\rangle = \langle \alpha, x\rangle$ for all $x \in \mathrm{postLie}(C)$.*

The GL-exponential is similarly defined from the GL product $*$:

**Definition 20 (GL-exponential).** The *GL-exponential* $\exp^* \colon \mathfrak{g}_{\mathrm{LB}} \to G_{\mathrm{LB}}$ is defined as

$$\exp^*(\alpha) = \mathbb{I} + \alpha + \frac{1}{2}\alpha * \alpha + \frac{1}{6}\alpha * \alpha * \alpha + \cdots = \sum_{j=0}^{\infty} \frac{1}{j!}\alpha^{*j}. \tag{81}$$

Recursive formulas for the coefficients of $\exp^*(\bullet)$ are found in [25, 19]. Here we derive a remarkably simple recursion formula based on the magmatic decomposition of OF, to our knowledge not found elsewhere:

**Lemma 3.** *For $\omega = \omega_L \times_\bullet \omega_R$ we have*

$$\langle \exp^*(\bullet), \mathbb{I}\rangle = 1 \tag{82}$$

$$\langle \exp^*(\bullet), \omega\rangle = \frac{1}{|\omega|} \cdot \langle \exp^*(\bullet), \omega_L\rangle \cdot \langle \exp^*(\bullet), \omega_R\rangle, \tag{83}$$

*or equivalently*

$$\langle \exp^*(\bullet), \omega\rangle = \frac{1}{\omega¡}, \tag{84}$$

*where $\omega¡$ denotes the ordered forest exponential.*

*Proof.* The derivation $D\exp^*(\bullet)$ satisfies $\langle D\exp^*(\bullet), \omega\rangle = |\omega|\langle \exp^*(\bullet), \omega\rangle$. On the other hand, since the $t$-scaling of the exponential is $t \star \exp^*(\bullet) = \exp^*(t\bullet)$ we find

$$D\exp^*(\bullet) = \left.\frac{d}{dt}\right|_{t=1} \exp^*(t\bullet) = \exp^*(t\bullet) * \bullet|_{t=1} = \exp^*(\bullet) * \bullet = \exp^*(\bullet)(\exp^*(\bullet) \rhd \bullet),$$

where we in the rightmost equality use (32) and $\Delta_{\sqcup\!\sqcup}(\exp^*(\bullet)) = \exp^*(\bullet) \otimes \exp^*(\bullet)$, since $\exp^*(\bullet) \in G_{\mathrm{LB}}$. Since $\omega_L \times_\bullet \omega_R = \omega_L(\omega_R \rhd \bullet)$ we find

$$\begin{aligned}\langle \exp^*(\bullet), \omega\rangle &= \frac{1}{|\omega|} \cdot \langle D\exp^*(\bullet), \omega\rangle = \frac{1}{|\omega|} \cdot \langle \exp^*(\bullet)(\exp^*(\bullet) \rhd \bullet), \omega_L(\omega_R \rhd \bullet)\rangle \\ &= \frac{1}{|\omega|} \cdot \langle \exp^*(\bullet), \omega_L\rangle \cdot \langle \exp^*(\bullet), \omega_R\rangle.\end{aligned}$$

□

The exponential is thus given as

$$\exp^*(\bullet) = \sum_{\omega \in \mathrm{OF}} \frac{\omega}{\omega¡}, \tag{85}$$

which justifies the naming of ¡ as a factorial function.

The computation of $\exp^*(\alpha)$ for an arbitrary $\alpha \in \mathfrak{g}_{\mathrm{LB}}$ can be done by the susbstitution: If $a(\bullet) = \alpha$ then

$$\begin{aligned}\langle \exp^* \alpha, \omega\rangle &= \langle \exp^* a(\bullet), \omega\rangle = \langle \exp^*(a \star \bullet), \omega\rangle \\ &= \langle a \star \exp^*(\bullet), \omega\rangle = \langle \exp^*(\bullet), a^t_\star(\omega)\rangle = \frac{1}{a^t_\star(\omega)¡},\end{aligned}$$

where the forest exponential ¡ is extended to polynomials by linearity.

**Backward error.**

Whereas $\exp^* \colon \mathfrak{g}_{\mathrm{LB}} \to G_{\mathrm{LB}}$ computes the exact flow operator, the inverse $\log^* \colon G_{\mathrm{LB}} \to \mathfrak{g}_{\mathrm{LB}}$ inputs a flow map, and computes the vector field generating this flow. In numerical analysis this is called the backward error analysis operator and is an important tool for analysing numerical integrators. The GL-logarithm $\log^*$ is defined for $\alpha \in G_{\mathrm{LB}}$ as

$$\log^* \alpha = \sum_{n=1}^{\infty} \frac{(-1)^{n-1}}{n} (\alpha - \delta)^{*n},$$

where $\delta \in G_{\mathrm{LB}}$ is the identity in the Lie–Butcher group, given as $\langle \delta, \mathbb{I}\rangle = 1$ and $\langle \delta, \omega\rangle = 0$ for $\omega \in \mathrm{OF}_C \backslash \{\mathbb{I}\}$. The GL-logarithm can be computed via its dual operation, the *eulerian idempotent* $e \in \mathrm{End}(\mathrm{k}\langle \mathrm{OF}_C\rangle)$ such that

$$\langle \log^*(\alpha), \omega\rangle = \langle \alpha, e(\omega)\rangle.$$

To compute $e$, we introduce the augmented GL-coproduct defined as

$$\overline{\Delta}_*(\omega) := \Delta_*(\omega) - \omega \otimes \mathbb{I} - \mathbb{I} \otimes \omega.$$

The recursion for $\Delta_*(\omega)$ (66)-(67) yields the following recursion for $\overline{\Delta}_*(\omega)$:

$$\overline{\Delta}_*(\mathbb{I}) = -\mathbb{I} \otimes \mathbb{I} \tag{86}$$

$$\overline{\Delta}_*(\omega_L \times_c \omega_R) = (\overline{\Delta}_*(\omega_L) + \omega_L \otimes \mathbb{I}) \sqcup\!\sqcup \times_c (\overline{\Delta}_*(\omega_R) + \omega_R \otimes \mathbb{I}). \tag{87}$$

The eulerian idempotent is computed as

$$e(\omega) = \sum_{n \geq 1} \frac{(-1)^{n-1}}{n} \sqcup\!\sqcup_n \overline{\Delta}_*^{n-1}(\omega),$$

where $\sqcup\!\sqcup_n$ is the shuffle of $n$ arguments and $\overline{\Delta}_*^n$ is the $n$-fold repeated application of the augmented GL coproduct. See Table 7 on p.41 for calculations of the eulerian idempotent for all forests up to and including order 4.

Since $\alpha$ is a character, we obtain the following formula for the backward error

$$\langle \log^*(\alpha),\omega\rangle = \sum_{n\geq 1}\frac{(-1)^{n-1}}{n}\sum_{\overline{\Delta}_*^{n-1}(\omega)}\langle \alpha,\omega_{(1)}\rangle\cdot\langle\alpha,\omega_{(2)}\rangle\cdots\langle\alpha,\omega_{(n)}\rangle. \tag{88}$$

**The development.**

For a curve $y(t)$ on a Lie group $G$, the *development* is a curve $\gamma(t)\in\mathfrak{g}$ such that $y'(t)=\gamma(t)y(t)$, thus $\gamma(t)=y'(t)y(t)^{-1}$ is given by the logarithmic derivative. There is a corresponding[9] combinatorial operation on $G_{\mathrm{LB}}$, given by a linear map $L\colon \mathrm{k}\langle \mathrm{OT}_C\rangle\to\mathrm{k}\langle\mathrm{OT}_C\rangle$ called the *Dynkin operator*, such that

$$\langle \alpha^{\cdot -1}\cdot D\alpha,\omega\rangle = \langle\alpha, L(\omega)\rangle \qquad \text{for every } \alpha\in G_{\mathrm{LB}}. \tag{89}$$

**Lemma 4.** *The Dynkin operator $L$ is computed as a convolution of endomorphisms, $L, S_\cdot, D\in\mathrm{End}(\mathscr{H}'_\cdot)$,*

$$L = S_\cdot * D := ш(S_\cdot\otimes D)\Delta_\cdot,$$

*where $\mathscr{H}'_\cdot$ is the Hopf algebra on $\mathrm{k}\langle\mathrm{OT}_C\rangle$ with shuffle ш as product, de-concatenation $\Delta_\cdot$ coproduct and antipode $S_\cdot$, and with grading $|\omega|$ counting nodes in the forest. Explicitly we have*

$$L(\omega)=\sum_{\Delta_\cdot(\omega)} S_\cdot(\omega_{(1)})\,ш\,\omega_{(2)}|\omega_{(2)}|. \tag{90}$$

*Proof.*

$$\begin{aligned}\langle\alpha^{\cdot-1}\cdot D\alpha,\omega\rangle &= \langle\alpha^{\cdot-1}\otimes D\alpha,\Delta_\cdot\omega\rangle = \sum_{\Delta_\cdot(\omega)}\langle\alpha,S_\cdot(\omega_{(1)})\rangle\langle\alpha,D(\omega_{(2)})\rangle\\ &= \langle\alpha, S_\cdot(\omega_{(1)})\,ш\,D\omega_{(2)}\rangle = \langle\alpha,(S_\cdot * D)(\omega)\rangle.\end{aligned}$$

□

Table 7 on p.41 contain the Dynkin map applied to all ordered forests up to and including order 4.

The inverse of the Dynkin map, denoted $\mathrm{evol}\colon\mathfrak{g}_{\mathrm{LB}}\to G_{\mathrm{LB}}$, yields a formal LB-series solution to equations of Lie type, $y'(t)=\gamma(t)y(t)$, for $y(t)\in G$, where $\gamma(t)\in\mathfrak{g}$ is given by a LB-series. In [11] it is proven that

$$\mathrm{evol}(\alpha)=\mathbb{I}+\sum_{n\geq 1}\sum_{\substack{n_1+\cdots+n_k=n\\ n_j>0}}\frac{\alpha_{n_1}*\alpha_{n_2}*\cdots*\alpha_{n_k}}{n_1(n_1+n_2)\cdots(n_1+n_2+\cdots+n_k)},$$

where $\alpha=\sum_{k\geq1}\alpha_k$ and $|\alpha_k|=k$ and $*$ is the convolution in $\mathscr{H}'_\cdot$. For $\omega\in\mathrm{OF}_C\backslash\{\mathbb{I}\}$ this yields

[9] Since the action of differentiation operators composes contravariantly, the order of right and left is swapped in the mapping from LB-series to differential equations on manifolds.

$$\langle \mathrm{evol}(\alpha), \omega \rangle = \sum_{n\geq 1} \sum_{\Delta_\cdot^{n-1}(\omega)} \frac{\langle \alpha, \omega_{(1)} \rangle \cdot \langle \alpha, \omega_{(2)} \rangle \cdots \langle \alpha, \omega_{(n)} \rangle}{|\omega_{(1)}| \cdot \big(|\omega_{(1)}| + |\omega_{(2)}|\big) \cdots \big(|\omega_{(1)}| + |\omega_{(2)}| + \cdots + |\omega_{(n)}|\big)},$$

and from this we find the recursion formulae

$$\langle \mathrm{evol}(\alpha), \mathbb{I} \rangle = 1 \tag{91}$$

$$\langle \mathrm{evol}(\alpha), \omega \rangle = \frac{1}{|\omega|} \sum_{\Delta_\cdot(\omega)} \langle \mathrm{evol}(\alpha), \omega_{(1)} \rangle \cdot \langle \alpha, \omega_{(2)} \rangle \quad \text{for } \omega \in \mathrm{OF}_C \backslash \{\mathbb{I}\}. \tag{92}$$

# 5 Concluding remarks

In this paper we have summarized the algebraic structures behind Lie–Butcher series. For the purpose of computer implementations, we have derived recursive formulae for all the basic operations on Lie–Butcher series that have appeared in the literature over the last decade. The simplicity of the recursive formulae are surprising to us. The GL-coproduct, the GL-exponential, the backward error and the inverse Dynkin map are in our opinion significantly simpler in their recursive formulations than the direct.

## *5.1 Programming in Haskell*

We are in the process of making a software library for computations with post-Lie algebras and Lie-Butcher series. As we have seen in this paper, many of the structures and operations have nice recursive definitions. Functional programming languages are well suited for this type of implementation. Haskell is one of the most popular functional programming languages, it is named after the logician Haskell B. Curry. The development of Haskell started in 1987 after a meeting at the conference on *Functional Programming Languages and Computer Architecture* (FPCA 87), where the need for common language for research in functional programming languages was recognized. Haskell has since grown into a mature programming language, not only used in functional programming research but also in the industry.

Not only do Haskell encourage recursive definitions of functions, it also has algebraic data types which give us the opportunity to define recursive data types.

Functional programming language will usually result in shorter and more precise code compared to imperative languages. Mathematical ideas are often straightforward to translate into a functional language.

A feature of Haskell that come in handy when working with infinite structures is lazy evaluation, meaning that an expression will not be computed before it is needed. This is an excellent feature for working with Lie-Butcher series, since these are infinite series. The infinite series can only be evaluated on finite data, and when such a computation is requested the system performs the necessary intermediate computations.

Mathematical ideas such as functors and monads are very important concept in Haskell, for example IO in Haskell is implemented as a monad. Another example is the vector space constructor in Haskell is a monad, which makes it very easy to linear extend a function on basis element to a linear function between vector spaces. Two other examples of monads are the free functor and the universal enveloping functor. The elementary differential map of B-series and Lie–Butcher series fits also nicely into this picture.

Finally, we remark that the proof assistant Coq can output Haskell code, so for critical parts of the software one can prove correctness of the implementation in Coq and then output this as verified Haskell code.

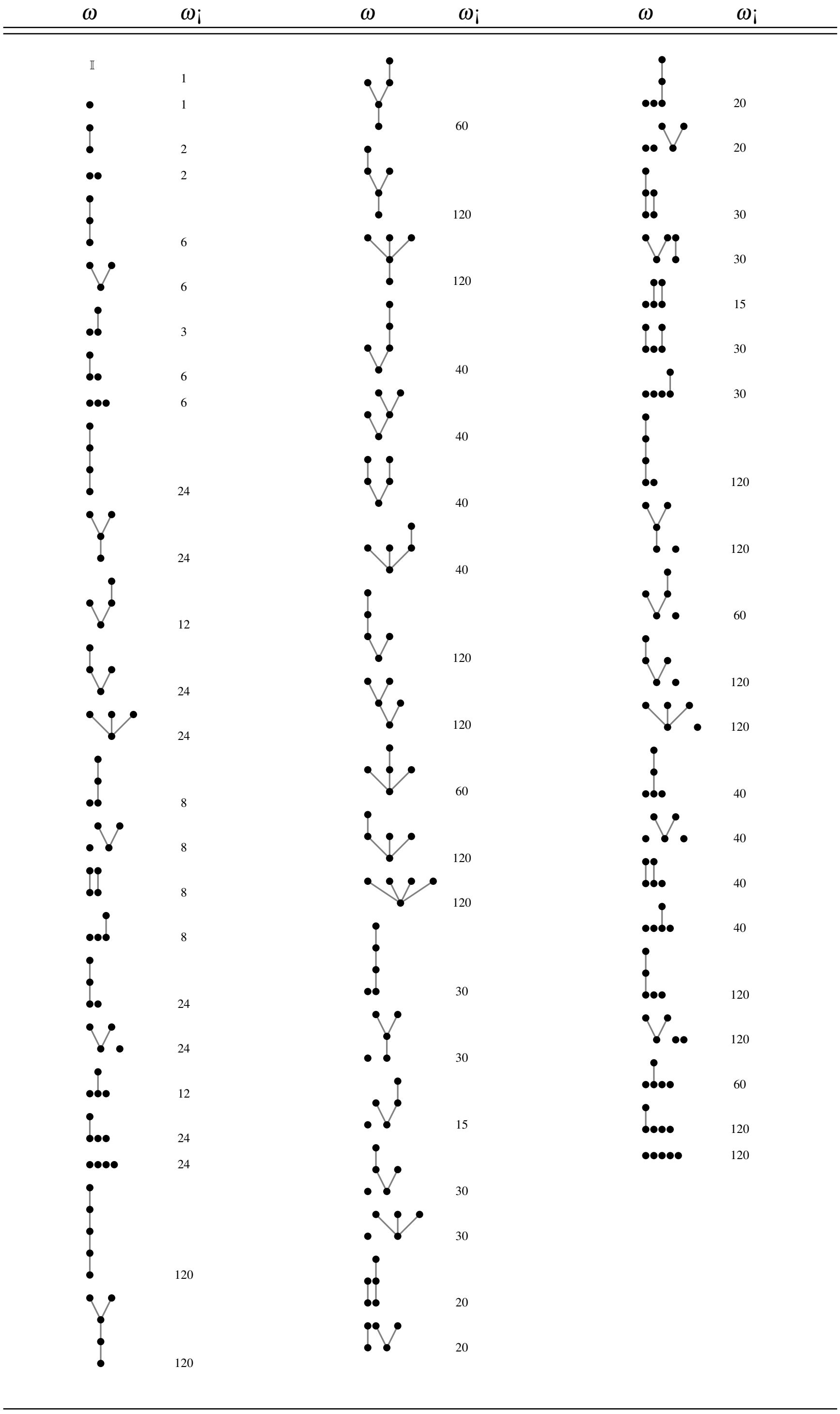

Table 1: Ordered forest factorial for all forest up to and including order 5.

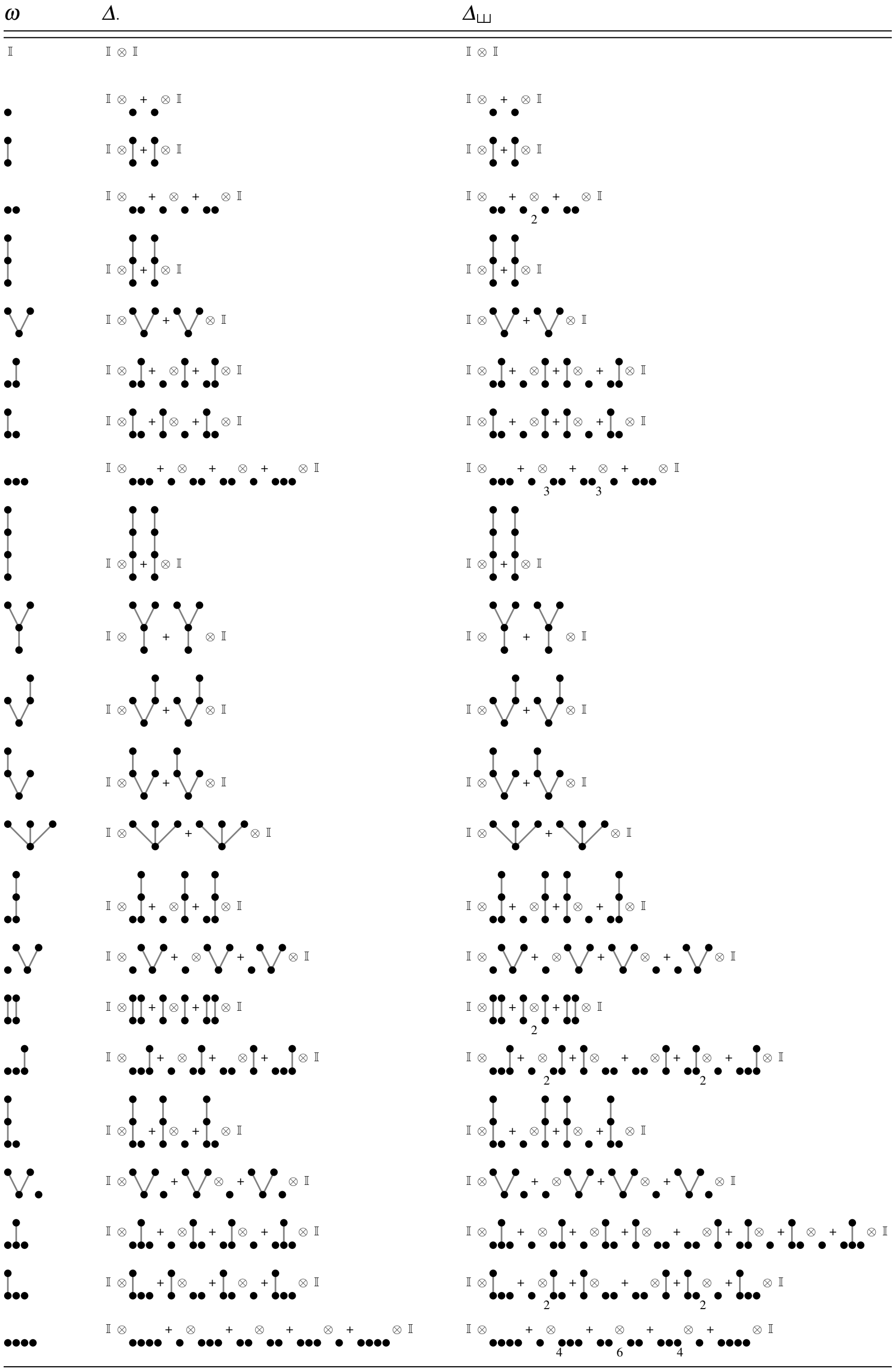

Table 2: Deconcatenation and deshuffle for ordered forest up to and including order 4. Note that the numbers under the terms are the coefficients to the terms.

| $\omega_1 \otimes \omega_2$ | $\omega_1 \rhd \omega_2$ | $\omega_1 * \omega_2$ |
|---|---|---|

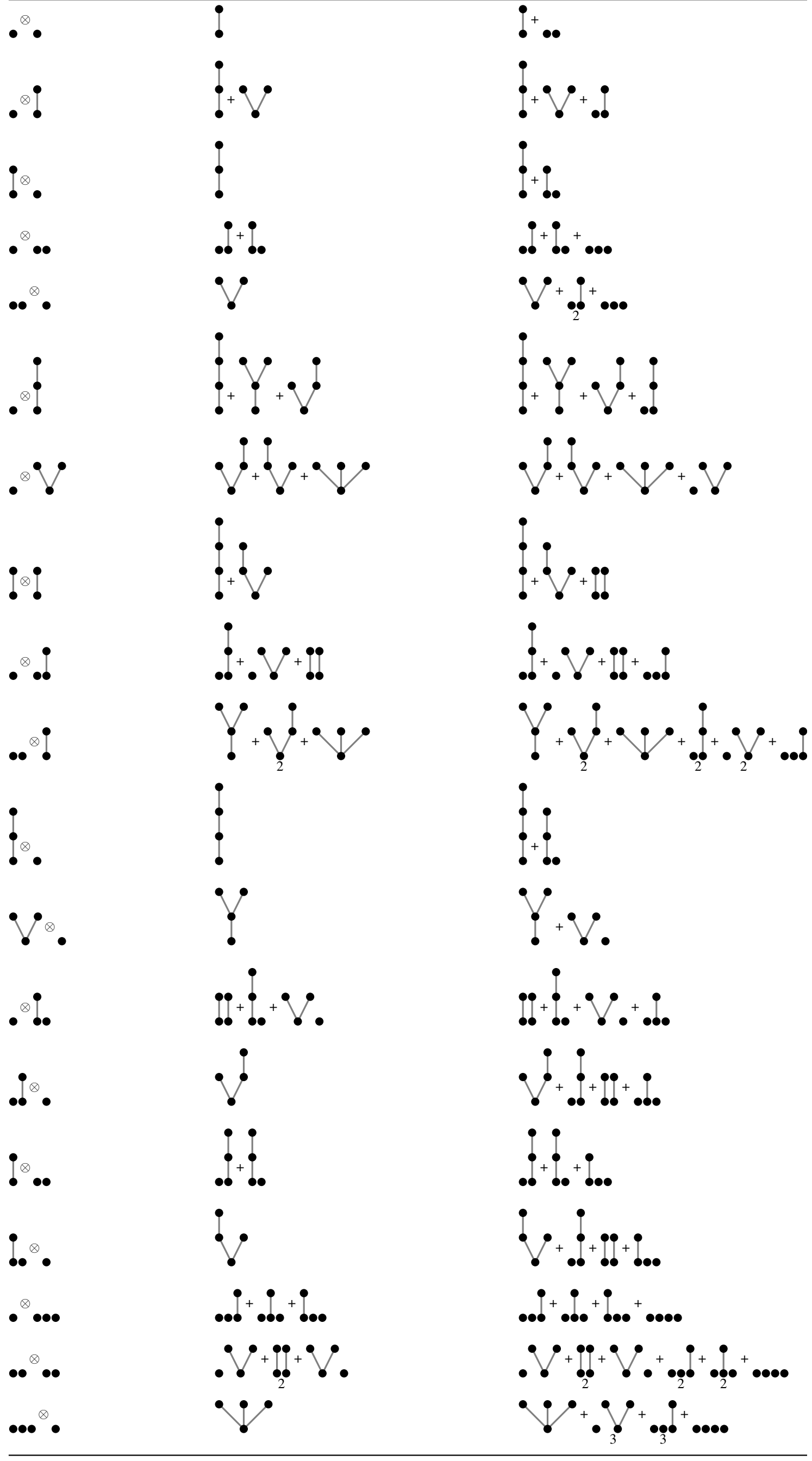

Table 3: Grafting and Grossman-Larsson product for all combinations of non-empty trees with total order up to and including order 4. Note that the numbers under the terms are the coefficients to the terms.

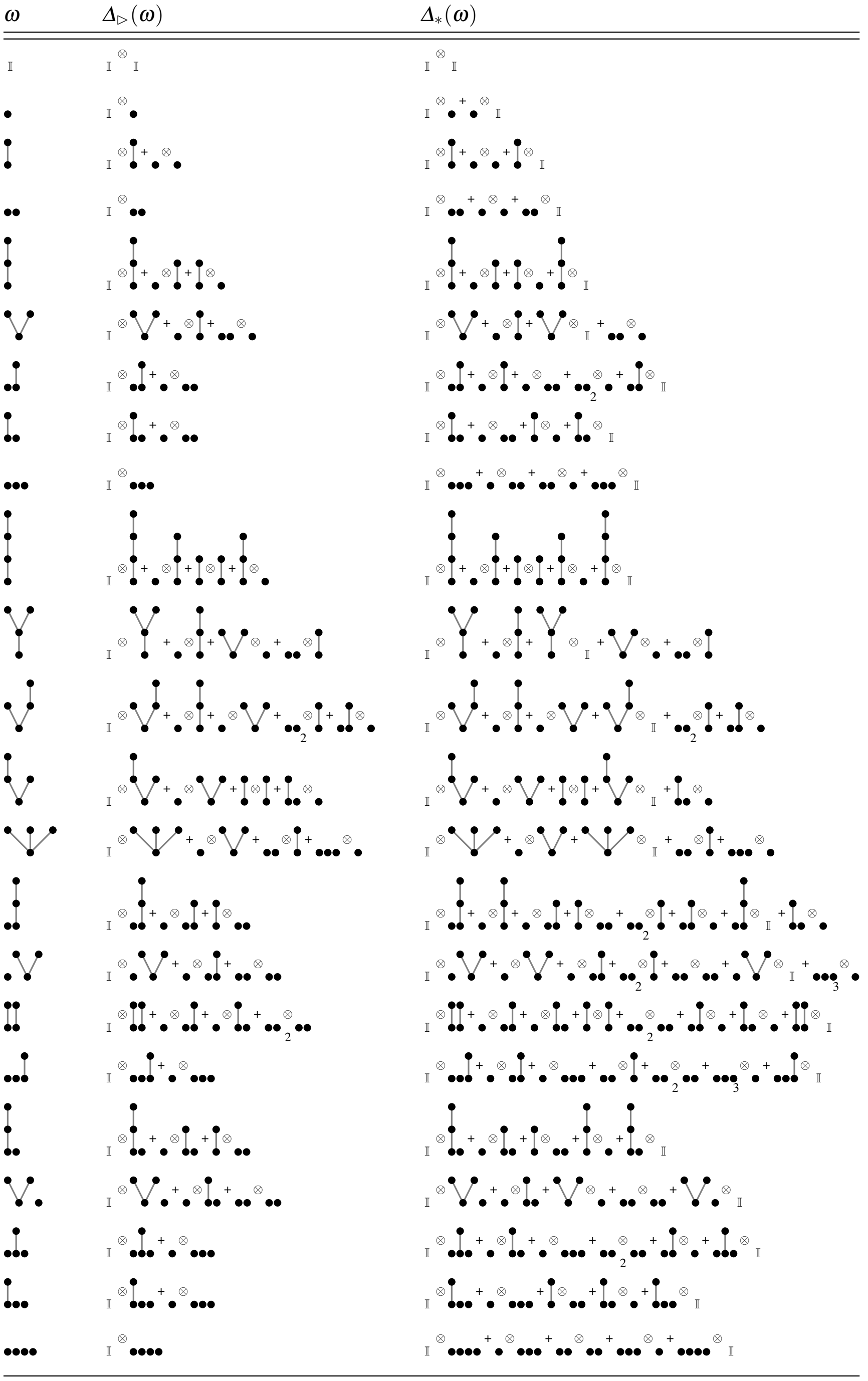

Table 4: Pruning and dual Grossman-Larsson coproduct for all forests up to and including order 4. Note that the numbers under the terms are the coefficients to the terms.

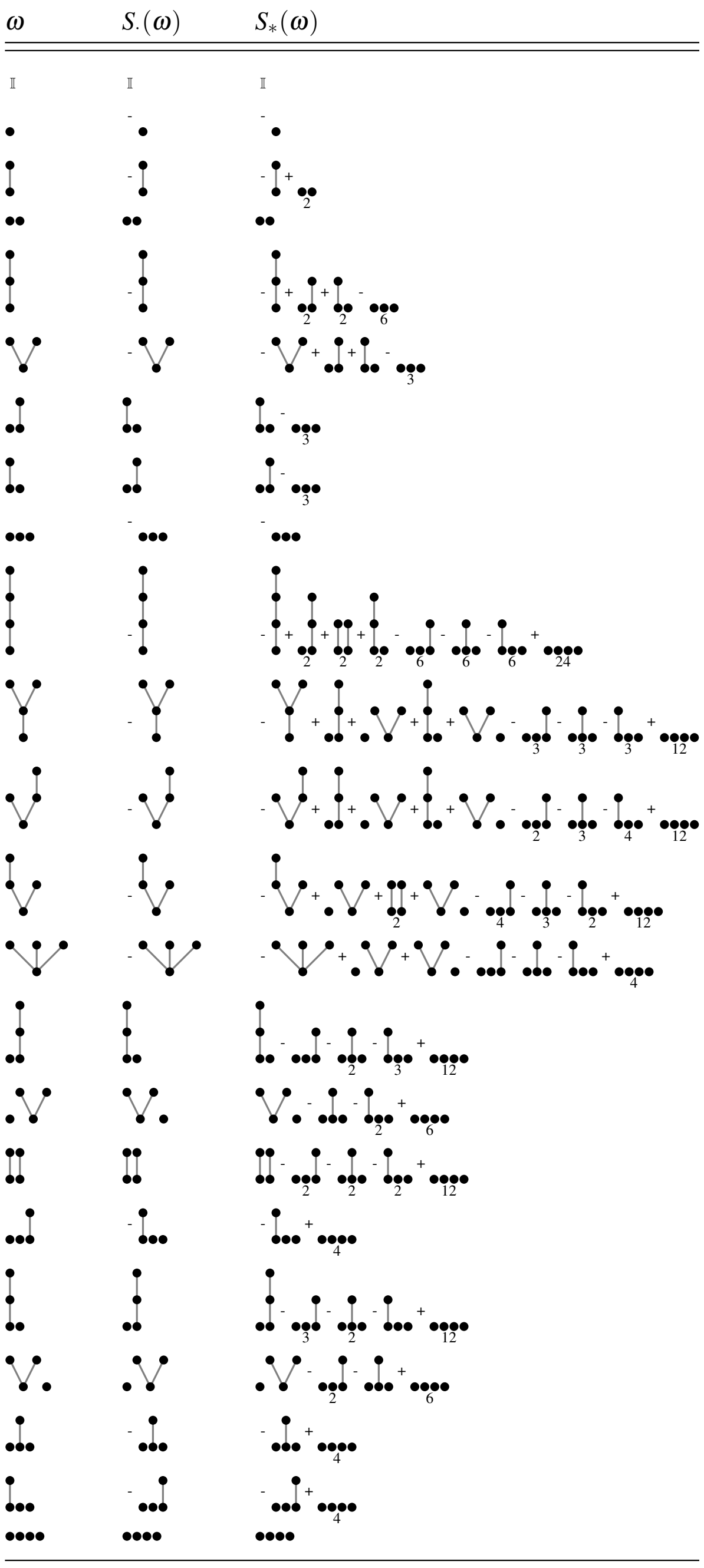

Table 5: Concatenation and Grossman-Larsson antipode map for all forests up to and including order 4. Note that the numbers under the terms are the coefficients to the terms.

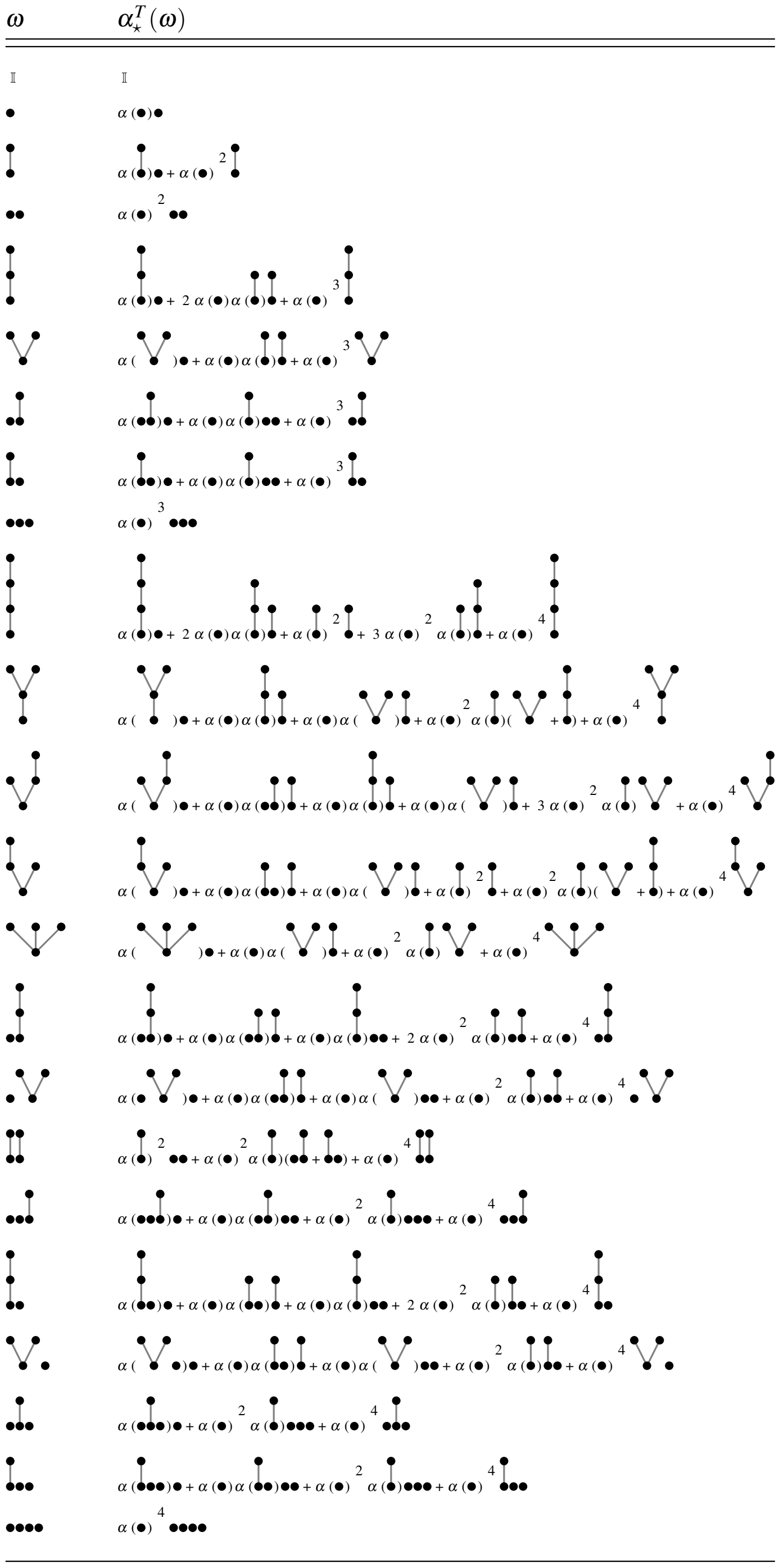

Table 6: Cosubstitution for an infinitesimal character $\alpha$ for all forest up to and including order 4.

| $\omega$ | $L(\omega)$ | $e(\omega)$ |
| --- | --- | --- |

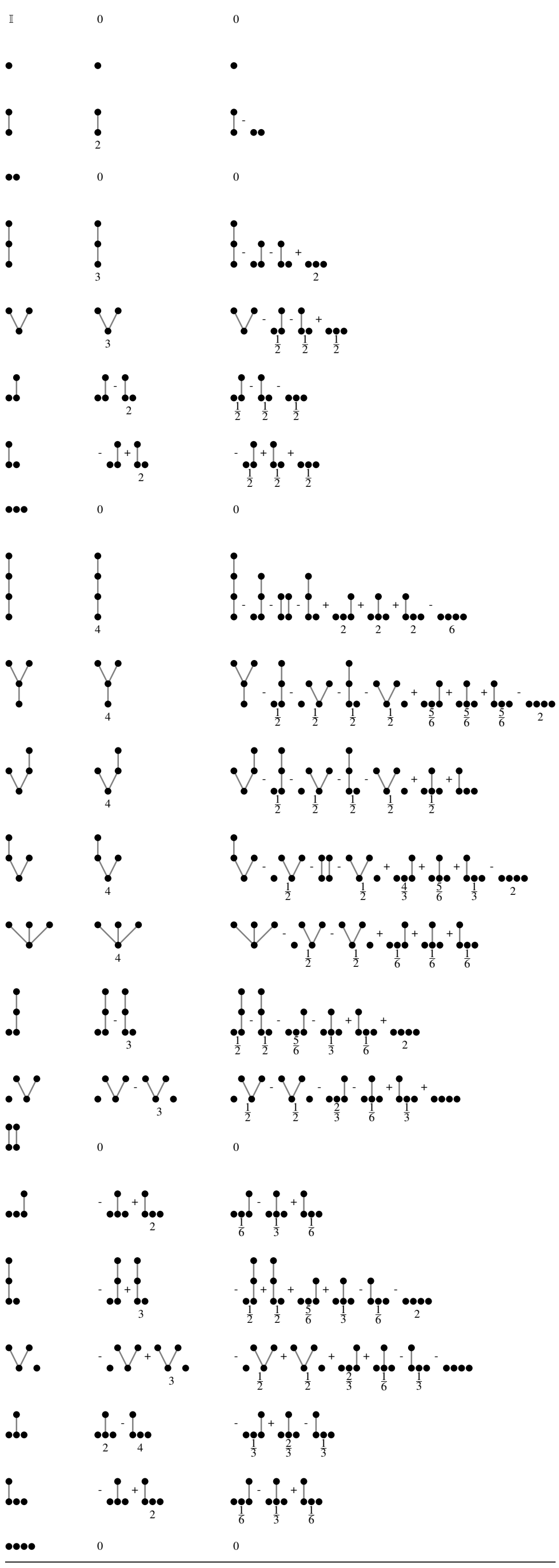

Table 7: Dynkin map $L$ and Eulerian idempotent $e$ for all forest up to and including order 4. Note that the numbers under the terms are the coefficients to the terms.